\documentclass[a4paper,11pt,twoside]{article}
\usepackage{pst-fill,pst-grad}
\usepackage{textcomp}
\usepackage[english]{babel}
\usepackage[utf8x]{inputenc}
\usepackage{graphicx}
\usepackage{amsmath}
\usepackage{multicol}
\usepackage{float}
\usepackage{fancyhdr,caption}
\usepackage[matrix,arrow,curve]{xy}
\usepackage{pstricks} 
\usepackage{amsmath,amsfonts,verbatim,afterpage,theorem,euscript,mathrsfs,amssymb}
\usepackage{amsfonts}
\usepackage{amssymb}
\usepackage{array}
\usepackage{dsfont}
\usepackage{hyperref}
\newcommand \A[1]{{\bf (#1)}}

\newcommand{\mysection}{\setcounter{equation}{0} \section}

\newtheorem{Definition}{Definition}[section]
\newtheorem{Proposition}{Proposition}[section]
\newtheorem{Lemme}{Lemma}[section]
\newtheorem{Theoreme}{Theorem}
\newtheorem{Corollaire}{Corollary}[section]
\newtheorem{Remarque}{Remark}[section]
\newcommand\R{\mathbb{R}}
\title{\bf Non Linear Singular Drifts and Fractional Operators}
\author{Diego Chamorro\footnote{Laboratoire de Math\'ematiques et Mod\'elisation d'Evry, Universit\'e d'Evry Val d'Essonne, France.}, St\'ephane Menozzi\footnote{Laboratoire de Math\'ematiques et Mod\'elisation d'Evry, Universit\'e d'Evry Val d'Essonne, France \& Laboratory of Stochastic Analysis, HSE, Moscow, Russia.}} 
\begin{document}
\maketitle
\begin{scriptsize}
\abstract{We consider parabolic PDEs associated with fractional type operators drifted by non-linear singular first order terms. When the drift enjoys some boundedness properties in appropriate Lebesgue and Besov spaces, we establish by exploiting a priori Besov-type estimates, the H\"older continuity of the solutions. In particular, we handle the \textit{almost} critical case in whole generality}.\\

\end{scriptsize}
\mysection{Introduction}

For a function $\theta:[0,+\infty[\times \mathbb{R}^{n}\longrightarrow\mathbb{R}$ with $n\geq 2$, we consider the following equation:
\begin{equation}\label{EquationPrincipaleIntro}
\begin{cases}
\partial_t \theta-\nabla\cdot(\mathbb{A}_{[\theta]}\,\theta)+\mathcal{L}^{\alpha}\theta=0,\qquad div(\mathbb{A}_{[\theta]})=0, \qquad 0<\alpha<2,\\[3mm]
\theta(0,x)=\theta_0(x), \quad x\in \R^n,
\end{cases}
\end{equation}
where $\mathcal{L}^{\alpha}$ is an $\alpha$-stable diffusion operator and where $\mathbb{A}_{[\theta]}=\mathbb{K}(\theta)$ is a \emph{nonlinear divergence free} velocity field  which is given as a vector $\mathbb{K}$ of singular integrals of convolution type in the space variable (see Section \ref{Secc_HipoetTheoremes} below for a precise definition of these objects). Observe that the general singular integral setting allows to consider a wide range of non-linear drifts with remarkable properties studied in many books (see e.g. \cite{Meyer, Stein2, Torchinski}). Also, the divergence free condition for the drift is very natural in problems arising from fluid dynamics. \\

When $\mathcal{L}^{\alpha}=(-\Delta)^{\frac{\alpha}{2}}$ we mention as a key example the surface quasi-geostrophic equation (SQG) where $n=2$ and for which we have
$$\mathbb{K}(\theta)=[-R_2\theta, R_1\theta],$$ 
where $(R_j)_{1\leq j\leq 2}$ denotes the Riesz transforms given by $\widehat {R_j \theta}= -i \frac{\xi_j}{|\xi|}\widehat \theta$. This equation has been considered from several points of view by many authors, see for example \cite{Caffarelli, CCCGW,PGDCH,CCZV,CW,Cordoba,Marchand,KN,Silv}, where problems of existence, uniqueness and regularity have been studied. Note that other examples of the previous setting derived from magnetohydrodynamic equations (MHD) have also been studied, see \cite{frie:vico:12}.\\

It is usual to decompose the study of this type of equations following the values of the fractional power of the Laplacian (or more generally following the smoothness degree $\alpha$ of the operator $\mathcal{L}^{\alpha}$). Indeed, if $1<\alpha<2$ then the regularizing effect of the operator $(-\Delta)^{\frac{\alpha}{2}}$ is strong enough to consider singular drifts for existence and regularity issues and this case is know as \emph{sub-critical}, for example singular drifts in Morrey spaces were considered in \cite{Wu}. The \emph{critical} case is given when $\alpha=1$, here the regularizing effect of the diffusion operator coincides with the derivative of the drift and H\"older regularity properties of the solutions were obtained in \cite{KN} with a drift in the $BMO$ space. Finally, the \emph{super-critical} case is given when $0<\alpha<1$ and some smoothness on the drift is asked in order to obtain some H\"older regularity for the solutions, see \cite{CW}. \textcolor{black}{We can refer as well to the works \cite{silv:12}, \cite{chau:meno:prio:20} for linear equations which clearly emphasize how in the supercritical case the smoothness of the drift is needed to derive a Schauder type theory}.\\

This division into cases is mainly driven by \emph{homogeneity} arguments which give a natural framework to study this type of equations, and this is valid for \emph{linear} drifts as well as for \emph{nonlinear} drifts. Indeed, in the case of \emph{linear} equations (when the drift $\mathbb{A}$ is independent from the solution):
\begin{equation}\label{EquationLineaireIntro}
\begin{cases}
\partial_t \psi\pm \nabla\cdot(\mathbb{A}\,\psi)+\mathcal{L}^{\alpha}\psi=0,\qquad div(\mathbb{A})=0,\qquad0<\alpha<2,\\[3mm]
\psi(0,x)=\psi_0(x),\qquad x\in \R^n,
\end{cases}
\end{equation}
\quad\\
we proved in \cite{DCHSM} that if the drift $\mathbb{A}$ is bounded (in the space variable) in Morrey-Campanato\footnote{See \cite{Adams} and \cite{Adams2} for a precise definition and for more properties of Morrey-Campanato spaces.} spaces $M^{q,a}(\R^{n})$ where the parameters $q,a$ are related to the dimension $n$ and to the smoothness degree $\alpha$ of the diffusion operator $\mathcal{L}^{\alpha}$ by the expression
\begin{equation}\label{ConditionEqLineaire}
\frac{a−n}{q}=1−\alpha,
\end{equation}
then from an initial data $\psi_0\in L^{\infty}(\R^{n})$ it is possible to derive some H\"older regularity $\mathcal{C}^{\gamma}(\R^{n})$ for the solutions of (\ref{EquationLineaireIntro}) where $0<\gamma<\alpha$. Remark in particular that, following this relationship above, if $1<\alpha<2$ then the corresponding space $M^{q,a}(\R^{n})$ can contain singular objects (recall that $L^{q}(\R^{n})\subset M^{q,0}(\R^{n})$), while if $\alpha=1$ we should have $a=n$ and then the corresponding space is $M^{q,n}(\R^{n})\simeq BMO(\R^{n})$. Finally, if $0<\alpha<1$, we have the identification $M^{q,a}(\R^{n})\simeq \mathcal{C}^{1-\alpha}(\R^{n})$, which corresponds with the smoothness asked for the drift in order to obtain a gain of regularity for the solutions. Thus, in all the cases stated above (sub-critical, critical and super-critical), as long as the drift is bounded in Morrey-Campanato spaces and satisfies the relationship (\ref{ConditionEqLineaire}) then it is possible to obtain a gain of regularity for the solutions. However, if this condition is not fulfilled then counterexamples to this gain of regularity can be produced, see \cite{Silv}.\\

As we can see, when looking for a gain of regularity, the linear case is very rigid as it is not possible to by-pass the condition (\ref{ConditionEqLineaire}). However, if we consider \emph{nonlinear} equations as (\ref{EquationPrincipaleIntro}) where $\mathbb{A}=\mathbb{A}_{[\theta]}$ then, by exploiting suitable information, there is a hope to break down this relationship  (\ref{ConditionEqLineaire}) in the following sense:\\

\emph{For an initial data $\theta_{0}$ in some Lebesgue spaces (no regularity asked for the initial data) and for a nonlinear drift $\mathbb{A}_{[\theta]}$ that satisfies some boundedness hypotheses, obtain a \emph{gain} of regularity for the solutions $\theta(t,x)$ (in terms of H\"older spaces in the space variable) for a smoothness degree $\alpha$ \emph{smaller} than the one given by the ``homogeneity" condition (\ref{ConditionEqLineaire})}.\\

One particular example of this situation is given in the celebrated article \cite{Caffarelli} that studies the SQG equation with a initial data $\theta_{0}\in L^{2}(\R^2)\subset M^{2,0}(\R^2)$. Following the relationship (\ref{ConditionEqLineaire}) since $n=2$ we should have $\alpha=2$, but with a careful treatment of the nonlinear drift it is shown that $\alpha=1$ is enough to obtain a gain of regularity for the solutions. Remark that this was done in two steps: first there is a ``de-singularization'' procedure in the space variable from $L^\infty_tL_x^{2}$ to $L^\infty_tL_x^{\infty}$ and then a h\"olderian gain of regularity (in the space variable) is deduced for a corresponding smoothness degree $\alpha=1$.\\

Another example was discussed in \cite{DCHSM_POTA} in the sub-critical case $\frac54<\alpha<2$ where it was possible to ask a little less regularity on the drift than the pure Morrey-Campanato condition (\ref{ConditionEqLineaire}). The main idea of this article was based in the remark that equation (\ref{EquationPrincipaleIntro}) admits a \emph{maximum principle} in terms of Lebesgue spaces: 
\begin{equation}\label{MaximumEstimate}
\|\theta(t,\cdot)\|_{L^{p}}\leq \|\theta_{0}\|_{L^{p}},
\end{equation}
for $2\leq p<+\infty$, but it \emph{also} admits an a priori ``energy" estimate given in terms of Besov spaces that reads as follows: 
\begin{equation}\label{AprioriBesovIntro}
\|\theta(\cdot,\cdot)\|_{L^{p}_t(\dot{B}^{\frac{\alpha}{p},p}_{p,x})}\leq C\|\theta_{0}\|_{L^{p}}\qquad \mbox{for } 2\leq p<+\infty,
\end{equation}
which provides us with a small regularity information of order $\frac{\alpha}{p}$ for the solution and this information can also be recovered for the nonlinear drift: indeed, in the nonlinear setting of equation (\ref{EquationPrincipaleIntro}), by assuming some mild boundedness properties for the drift $\mathbb{A}_{[\theta]}$ in Lebesgue (in the $t$ variable) and Besov (in the $x$ variable) spaces of the form
$$\|\mathbb{A}_{[\theta]}(\cdot,\cdot)\|_{L^{p}_t(\dot{B}^{\frac{\alpha}{p},p}_{p,x})}\leq C\|\theta(\cdot,\cdot)\|_{L^{p}_t(\dot{B}^{\frac{\alpha}{p},p}_{p,x})},$$
and coupling this estimate with the previous Besov control (\ref{AprioriBesovIntro}) it was possible to \textit{break} the homogeneity relation (\ref{ConditionEqLineaire}) through an interplay between the Morrey-Campanato and Besov stability properties of the drift. However, in the quoted work  \cite{DCHSM_POTA} we could only handle the sub-critical case $\frac54<\alpha<2$ up to an additional mollification of the drift in the time variable.  Indeed, we missed in \cite{DCHSM_POTA} an $L^\infty$ in time and Besov in space control which is precisely established here (see Proposition \ref{PropoLpLinfty1} below). This amends the proofs in \cite{DCHSM_POTA}.\\

The main goal of this article is to break the condition (\ref{ConditionEqLineaire}), in the sense introduced above, for the nonlinear equation \eqref{EquationPrincipaleIntro} in the sub-critical case where $\alpha=1+\varepsilon$ with $\varepsilon>0$ meant to be small, improving in this way the results  in \cite{DCHSM_POTA} and considering the almost critical case for a large class of non-linear drifts.\\

The underlying idea is to use the sub-criticality, \emph{i.e.} the fact that $\alpha=1+\varepsilon>1$, to deduce a uniform in time gain of regularity in terms of Sobolev or Besov spaces. Then, by using \emph{all} the a priori informations available (maximum principle and Besov uniform estimates) we can obtain a better control of the evolution of the profile of the solutions of the linear equation (\ref{EquationLineaireIntro}) from which we will deduce, by duality, a gain of regularity for the solutions $\theta(t,x)$ of the equation (\ref{EquationPrincipaleIntro}).
\mysection{Assumptions and Main Results}\label{Secc_HipoetTheoremes}

In order to state our theorems we need to be more precise about the objects that define the equation (\ref{EquationPrincipaleIntro}) as some properties of these objects are essential in this work.
\begin{itemize}
\item[\A{A}] {\bf The diffusion operator $\mathcal{L}^{\alpha}$.}
The L\'evy-type operator $\mathcal{L}^{\alpha}$ appearing in equation (\ref{EquationPrincipaleIntro}) can be viewed as a natural generalization of the usual fractional Laplacian $(-\Delta)^{\frac{\alpha}{2}}$. Indeed, for $0<\alpha<2$  and for a function $\varphi:\R^n\longrightarrow \R$ regular enough, we can define the operator $\mathcal{L}^{\alpha}$  in the Fourier variable by the expression
$$\widehat{\mathcal{L}^{\alpha}\varphi}(\xi)=a(\xi)\widehat{\varphi}(\xi),$$
where the symbol $a$ is of the form $\displaystyle{a(\xi)=\int_{\mathbb{R}^n\setminus \{0\}}\big(1-\cos(\xi\cdot y)\big)\pi(y)dy}$, 
with a symmetric function $\pi$, \emph{i.e.} $\pi(y)=\pi(-y)$, satisfying for all $y\in \R^n$ the following bounds:
\begin{equation}\label{DefKernelLevy}
\begin{split}
\mathfrak{c}_1|y|^{-n-\alpha}\leq &\pi(y)\leq \mathfrak{c}_2|y|^{-n-\alpha} \qquad \mbox{over } |y|\leq 1,\\[2mm]
0\leq &\pi(y)\leq \mathfrak{c}_2|y|^{-n-\alpha} \qquad \mbox{over } |y|> 1.
\end{split}
\end{equation}
Here $0<\mathfrak{c}_1\leq\mathfrak{c}_2$ are positive constants that are fixed once and for all. The important point here is that the parameter $\alpha$ (called the \emph{smoothness} or \emph{regularity} degree) will rule the smoothness properties of the operator $\mathcal{L}^{\alpha}$ in the sense that for such $\alpha$ the regularizing effect of $\mathcal{L}^{\alpha}$ is similar to the fractional Laplacian $(-\Delta)^{\frac{\alpha}{2}}$. See \cite{Jacob} for additional properties concerning Lévy operators. \\

\item[\A{B}] {\bf The nonlinear Drift $\mathbb{A}_{[\theta]}$.} First, for $1\leq i\leq n$, we set 
$$\displaystyle{\mathbb{K}_{i}(\theta)(t,x)=v.p. \int_{\mathbb{R}^{n}}\kappa_{i}(x-y)\theta(t,y)dy},$$
where $\kappa_{i}=\kappa_{i}(x)$ is the associated kernel of the singular integral operator $\mathbb{K}_{i}$ which acts only in the space variable: in the previous formula the dependence in the time variable only comes from the solution $\theta(t,\cdot)$. Then we consider the vector
\begin{equation}\label{DefinitionNoyau}
\mathbb{A}_{[\theta]}(t,x)=[\mathbb{K}_{1}(\theta),\cdots, \mathbb{K}_{n}(\theta)](t,x).
\end{equation}
\item[\A{C}] {\bf Boundedness conditions for the nonlinear drift term $\mathbb{A}_{[\theta]}$}. We set now the properties of the vector $\mathbb{A}_{[\theta]}$ given in (\ref{DefinitionNoyau}) that defines the nonlinear drift in equation (\ref{EquationPrincipaleIntro}) above. 
\begin{itemize}
 
\item[$\bullet$] \emph{Boundedness in Lebesgue spaces}: as only integrability in the space variable is studied here, for $1<p<+\infty$ we will need the following general condition
\begin{equation}\label{HypoLPBorne1}
\|\mathbb{A}_{[\theta]}(t, \cdot)\|_{L^p}\leq C\|\theta(t,\cdot)\|_{L^p}.
\end{equation}
Note that the condition (\ref{HypoLPBorne1}) amounts to ask for the space variable the boundedness assumption $\|\mathbb{K}_i(\theta)(t, \cdot)\|_{L^p}\leq C\|\theta(t, \cdot)\|_{L^p}$ for all the singular integral operators $\mathbb{K}_i$.\\

\item[$\bullet$] \emph{Boundedness in Besov spaces}: we study here the boundedness of the drift in terms of Besov spaces in the space variable. Thus for $1< p,q<+\infty$ and $0<\sigma<2$ we will assume the following boundedness property:
\begin{equation}\label{HypoBesovBorne01}
\|\mathbb{A}_{[\theta]}(t, \cdot)\|_{\dot{B}^{\sigma,p}_{q}}\leq C\|\theta(t,\cdot)\|_{\dot{B}^{\sigma,p}_{q}},
\end{equation}
and we will also assume the following permutation property:
\begin{equation}\label{HypoBesovBorne02}
\|(-\Delta)^{\frac{\sigma}{2}}\mathbb{A}_{[\theta]}(t, \cdot)\|_{L^p}=\|\mathbb{A}_{[(-\Delta)^{\frac{\sigma}{2}}\theta]}(t, \cdot)\|_{L^p},
\end{equation}
for all $1<p<+\infty$ and all $0<\sigma<2$.
\end{itemize}

It is worth noting here that all Calder\'on-Zygmund operators $\mathbb{K}_i$  (such as the Riesz transforms) satisfy these boundedness conditions and that \emph{no regularity is asked in the space variable} for $\mathbb{K}_i$. The positive constants above depend on the operator $\mathbb{A}$ and will be fixed once and for all.\\

\item[\A{D}] {\bf The initial condition $\theta_{0}$.} We will assume here that the initial data $\theta_{0}:\mathbb{R}^{n}\longrightarrow \mathbb{R}$ belongs to \emph{all} the Lebesgue space $L^{p}(\mathbb{R}^{n})$ with $1\leq p\leq {\bf\bar{p}}<+\infty$. Remark in particular that neither regularity nor boundedness are asked for the initial data. We will adopt from now on the following notation:
\begin{equation}\label{ControlNormeDonneeInitiale}
1\leq \mu=\underset{1\le p\leq {\bf\bar{p}}}{\max}\;\|\theta_{0}\|_{L^{p}}<+\infty,
\end{equation}
and the constant $\mu$ will be denoted as the \emph{size} of the initial data.
\end{itemize}
Within this framework we can state our results. We first present an existence result which is crucial for our analysis but somehow standard, see \emph{e.g.} \cite{DCHSM}, \cite{DCHSM_POTA}.
\begin{Theoreme}[The NonLinear case]\label{Theoreme1} Let $n\geq 2$. Under hypotheses \A{A}, \A{B}, \A{C} and \A{D}, for all fixed $T> 0$, there exists a  weak solution $\theta(t,x)$ to the nonlinear equation \eqref{EquationPrincipaleIntro} in $L^{\infty}([0,T], L^{p}(\mathbb{R}^n))\cap L^{p}([0,T], \dot{B}^{\frac{\alpha}{p},p}_p(\mathbb{R}^n))$ with ${\bf 2}\leq p\leq {\bf\bar{p}}$. Also, for all $0<t\leq T$ the maximum principle holds: 
\begin{equation}\label{MaximumPrinciple1}
\|\theta(t,\cdot)\|_{L^p}\leq \|\theta_0\|_{L^p}\leq \mu,
\end{equation}
and moreover we have the following Besov-norm a priori control
\begin{equation}\label{ControlBesovApriori1}
\|\theta\|_{L^{p}_t(\dot{B}^{\frac{\alpha}{p},p}_{p,x})}=\left(\int_{0}^{T}\|\theta(t,\cdot)\|_{\dot{B}^{\frac{\alpha}{p},p}_{p}}^{p}dt\right)^{\frac{1}{p}}\leq C \|\theta_0\|_{L^p}\leq C\mu.
\end{equation}
Furthermore, for $ \alpha=1+\varepsilon$, $0<\varepsilon \ll1$, if the index of integrability is big enough, say ${\bf \bar p}>\frac n\varepsilon$, there is a unique weak solution satisfying the above properties.
\end{Theoreme}
The existence of such solutions, as well as the associated controls \eqref{MaximumPrinciple1} and \eqref{ControlBesovApriori1}, can be derived proceeding similarly to Section 3 in \cite{DCHSM_POTA}. The point is to use compactness arguments through a vanishing viscosity approach combined with the stability conditions \eqref{HypoLPBorne1} and \eqref{HypoBesovBorne01}. Uniqueness is then derived in a second time, see Proposition \ref{Prop_UNIQUENESS_W_SOL} below, from the fact that any weak solution enjoying the above properties can be represented with an \textit{integral} formulation (see Lemma \ref{LEMME_MILD} below). This feature specifically comes from the fact that we are considering the sub-critical regime. Observe as well that, the integrability needed to derive uniqueness goes to infinity when approaching the critical case (\emph{i.e.} ${\bf \bar p}\to +\infty$ if $\varepsilon\to 0$).

Note that the restriction ${\bf 2}\leq p$ is mandatory to obtain the Besov-norm control (\ref{ControlBesovApriori1}) which is one of the key ingredients of this work. Indeed, this Besov-norm a priori information is no longer available in the case $1<p<2$, see \cite{PGDCH}. Remark also that the Besov regularity given in (\ref{ControlBesovApriori1}) is naturally linked to the smoothness degree $\alpha$ of the diffusion operator $\mathcal{L}^{\alpha}$ it is maximal and equal to $\frac{\alpha}{2}$ when $p=2$  and vanishes when $p\to +\infty$. \\

The main result of the current work is the following theorem.
\begin{Theoreme}[H\"older regularity]\label{THM_HOLDER}
Let the dimension $n\geq 2$ and consider over $\mathbb{R}^{n}$ the equation (\ref{EquationPrincipaleIntro}) where the diffusion operator $\mathcal{L}^{\alpha}$ satisfies \A{A} for some $\alpha$ such that $\alpha=1+\varepsilon$ with $0<\varepsilon \ll1$, the nonlinear drift $\mathbb{A}_{[\theta]}$ defined in \A{B} enjoys the stability properties of \A{C} and the initial data $\theta_{0}$ satisfies \A{D}. Consider any weak solution $\theta$ of \eqref{EquationPrincipaleIntro} in $\displaystyle{\bigcap_{{\bf 2}\leq p\leq{\bf \bar{p} }}}L^{\infty}([0,T], L^{p}(\mathbb{R}^n))\cap L^{p}([0,T], \dot{B}^{\frac{\alpha}{p},p}_p(\mathbb{R}^n))$, which for $\alpha=1+\varepsilon$ is unique if $ \bar {\bf p}\ge \frac{n}{\varepsilon}$, there exists a time $0<T_{0}<T$ such that for all $t>T_{0}$ the solution $\theta(t,\cdot) $ belongs to the H\"older space $\mathcal{C}^\gamma({\mathbb R}^n) $ for some $0<\gamma<1$.
\end{Theoreme}
As we can see, when $\alpha>1$ and although the initial data is not regular, we can deduce some H\"older regularity for weak solutions, however this gain of regularity is not instantaneous and some time is needed in order to obtain the h\"older continuity. The main arguments used here to prove this theorem are the following: first we need to prepare the information available on weak solutions, as at some point we will require to pass from an $L^p_t(\dot{B}^{\frac{\alpha}{p},p}_{p,x})$ control given by (\ref{ControlBesovApriori1}) to an $L^\infty_t(\dot{B}^{\frac{\sigma}{p},p}_{p,x})$ estimate. Again, the mild representation of Lemma \ref{LEMME_MILD} is here crucial. This first step is not absolutely trivial and the price to pay in order to obtain an $L^\infty$ control in the time variable can be seen in the space variable with a small loss of the regularity obtained (given by the fact that $\sigma<\alpha$), but this small loss of regularity can still be compensated by the regularity degree $\alpha>1$. Then with this $L^\infty$-Besov information at hand, we can study the evolution of weak solutions of the problem (\ref{EquationPrincipaleIntro}) via a dual equation (see Section \ref{Secc_EquationDuale}). The main argument relies then in a suitable control of the elements of the Hardy space which helps us characterize, by duality, the corresponding H\"older spaces.

Let us remark that when $\alpha\leq1$ the regularity loss observed in the first step in order to obtain the $L^\infty$-Besov information seems to be too strong and we believe that a further step is necessary in order to study the case $\alpha=1$ or the case $\alpha<1$. Anyhow, we feel that the approach below should work in the critical case $\alpha=1 $ should we start form a bounded initial condition (see also \emph{e.g.} the proof of Proposition \ref{Prop_UNIQUENESS_W_SOL} for related topics).
\mysection{Preliminar estimates for the nonlinear system}\label{SeccLpLinfty}
As already indicated above, a crucial tool for the whole analysis is the integral representation of weak solutions enjoying the properties indicated in Theorem \ref{Theoreme1}.
\begin{Lemme}[Integral representation of the weak solutions]\label{LEMME_MILD}
Let the dimension $n\geq 2$ and consider over $\mathbb{R}^{n}$ the equation (\ref{EquationPrincipaleIntro}) where the diffusion operator $\mathcal{L}^{\alpha}$ satisfies \A{A} with $\alpha=1+\varepsilon$, the nonlinear drift $\mathbb{A}_{[\theta]}$ defined in \A{B} enjoys the stability properties of \A{C} and the initial data $\theta_{0}$ satisfies \A{D}. Consider an associated weak solution $\theta(t,x)$ to the nonlinear equation \eqref{EquationPrincipaleIntro}. Then, the following integral representation of the solution holds. For all $t\ge 0 $,
\begin{equation}\label{REP_MILD_W_SOL}
\theta(t,x)=\mathfrak{p}^\alpha_t\ast \theta_0(x)+ \int_{0}^{t}\mathfrak{p}^\alpha_{t-s}\ast\nabla\cdot(\mathbb{A}_{[\theta]}\,\theta(s,x))ds,
\end{equation}
where $\mathfrak{p}^\alpha_t$ stands for the semi-group kernel associated with the operator $\mathcal{L}^\alpha $ at time $t$.
\end{Lemme}
The proof of this lemma, which mainly relies on properties of the underlying stable-type kernel,
 is postponed to Appendix \ref{ANNEX_MILD} to focus uniqueness and regularity issues. We anyhow provide now the proof of uniqueness of the weak solutions in the current setting which can be seen as a first important consequence of the previous \textit{integral} representation.
\subsection{On uniqueness of the weak solutions}
We have the following uniqueness result.
\begin{Proposition}\label{Prop_UNIQUENESS_W_SOL}
Under the 	assumptions of  Theorem \ref{Theoreme1}, for $\alpha=1+\varepsilon $ with $0<\varepsilon\ll 1$ and if ${\bf {\bar p}} > \frac n\varepsilon$, there exists a unique solution $\theta$ of \eqref{EquationPrincipaleIntro} in $\displaystyle{\bigcap_{{\bf 2}\leq p\leq{\bf \bar{p} }}}L^{\infty}([0,T], L^{p}(\mathbb{R}^n))\cap L^{p}([0,T], \dot{B}^{\frac{\alpha}{p},p}_p(\mathbb{R}^n))$.
\end{Proposition}
\textbf{Proof.} Let us start from two weak solutions $\theta_1,\theta_2$ to \eqref{EquationPrincipaleIntro} which belong to $L^\infty_t(L^p_x)\cap L_t^p(\dot B_{p,x}^{\frac \alpha p,p})$ with $p\in [2,{\bf \bar p}]$. Denoting by $\delta_{\theta}:=\theta_1-\theta_2$ and $\delta_{(\mathbb{A}_{[\theta]}\,\theta)}=\mathbb{A}_{[\theta_1]}\,\theta_1-\mathbb{A}_{[\theta_2]}\,\theta_2$  it thus holds from Lemma \ref{LEMME_MILD} above that for fixed $t$,
\begin{eqnarray*}
\|\delta_\theta(t,\cdot)\|_{\dot B_p^{\frac{\alpha}{p}, p}}&\le& \left\|\int_{0}^{t}(-\Delta)^{\frac{1}{2}}\mathfrak{p}^\alpha_{t-s}\ast  (\mathbb{A}_{[\theta_1]}\,\theta_1-\mathbb{A}_{[\theta_2]}\,\theta_2)(s,\cdot)ds\right\|_{\dot B_p^{\frac\alpha p,p}}\\
&\leq &\left\|\int_{0}^{t}(-\Delta)^{\frac{1}{2}}\mathfrak{p}^\alpha_{t-s}\ast\delta_{(\mathbb{A}_{[\theta]}\,\theta)}(s,\cdot)ds\right\|_{\dot B_{p}^{\frac \alpha p,p}}.
\end{eqnarray*}
Integrate next with respect to $t\in [0,T]$. We get:
\begin{eqnarray*}
\left(\int_0^{T }\|\delta_\theta(t,\cdot)\|_{\dot B_p^{\frac{\alpha}p, {p}}}^p dt\right)^{\frac 1p} &\le&\left\|\int_{0}^{\cdot}(-\Delta)^{\frac{1}{2}}\mathfrak{p}^\alpha_{\cdot-s}\ast \delta_{(\mathbb{A}_{[\theta]}\,\theta)}(s,\cdot)ds\right\|_{L_{t}^p(\dot B_{p,x}^{\frac \alpha p,p})}.
\end{eqnarray*}
Hence,
\begin{eqnarray*}
\left(\int_0^{T }\|\delta_\theta(t,\cdot)\|_{\dot B_{p}^{\frac{\alpha}{p}, {p}}}^p dt\right)^{\frac 1p} &&\le\left\|\int_{0}^{\cdot}(-\Delta)^{\frac{\alpha}{2}}\mathfrak{p}^\alpha_{\cdot-s}\ast (-\Delta)^{\frac{1-\alpha} 2}(\delta_{(\mathbb{A}_{[\theta]}\,\theta)}(s,\cdot))ds\right\|_{L_{t}^p(\dot B_{p,x}^{\frac \alpha p, p})}.
\end{eqnarray*}
Let us now recall the thermic characterization of the homogeneous Besov norm (see e.g. Triebel \cite{trie:83-2}):
\begin{equation}\label{CAR_THERM}
\|f\|_{\dot B_{p}^{\frac \alpha p,p}}^p\simeq \int_0^{+\infty}  \tau^{(1-\frac \alpha {2p})p}\|\partial_\tau \mathfrak h _\tau\ast f\|_{L^p(\R^n)}^p\frac {d\tau} \tau,
\end{equation} 
where $\mathfrak h $ stands for the standard (gaussian) heat kernel on $\R^n$.
 Introduce now 
\begin{equation*}
G^\alpha f(t,x)=\int_{0}^t P_{t-s}^\alpha f(s,x)ds,
\end{equation*}
the parabolic Green kernel, where $P^\alpha$ denotes the semi-group generated by $\mathcal L^\alpha $ (with density $\mathfrak p^\alpha $). Let us also define \begin{equation}
\label{DEF_PSI}
\Psi (s,\cdot )=(-\Delta) ^{\frac{1- \alpha} 2}(\delta_{(\mathbb A_{[\theta]}\theta)})(s,\cdot),
\end{equation}
Hence, using the Fubini theorem, we can apply  maximal-regularity type inequalities  to the kernel $(-\Delta)^{\frac{\alpha} 2}\mathfrak{p}^\alpha_{\cdot-s}$. Precisely we have,
\begin{align*}
\|\delta_\theta\|_{L_t^p(\dot B_{p,x}^{\frac \alpha p,p})}^p  \leq& C\int_0^{+\infty} \tau^{(1-\frac \alpha {2p})p} \int_0^T\left\|\partial_\tau \mathfrak h_\tau \ast (-\Delta)^{\frac{\alpha} 2}G^\alpha \Psi(t,\cdot)\right\|_{L^{p}}^p dt \frac {d\tau} \tau\notag\\
\le& C   \int_0^{+\infty} \tau^{(1-\frac \alpha {2p})p}\int_0^T  \|(-\Delta)^{\frac \alpha 2} G^\alpha \big( \Psi \ast \partial_\tau \mathfrak h_\tau\big)(t,\cdot)\|_{L^p}^pdt\frac {d\tau} \tau\notag\\
=&C   \int_0^{+\infty}  \tau^{(1-\frac \alpha {2p})p}\int_0^T  \|(-\Delta)^{\frac \alpha 2} G^\alpha \widetilde \Psi_\tau (t,\cdot)\|_{L^p}^p dt \frac {d\tau} \tau,
\end{align*}
with $\widetilde \Psi_\tau (t,\cdot):= \Psi \ast \partial_\tau \mathfrak h_\tau(t,\cdot)$ using the associativity of the convolution product. The maximal regularity now yields (recall that $p \in [2,{\bf \bar p}]$ with ${\bf {\bar p}} <+\infty$):
\begin{align*}
\|\delta_\theta\|_{L_t^p(\dot B_{p,x}^{\frac \alpha p,p})}^p  \leq& C\int_0^{+\infty}\tau^{(1-\frac \alpha {2p})p}\int_0^T  \|  \widetilde \Psi_\tau (t,\cdot)\|_{L^p}^p dt \frac {d\tau} \tau \notag\\
\leq& C\int_0^{+\infty} \tau^{(1-\frac \alpha {2p})p}\int_0^T  \|  \partial_\tau \mathfrak h_\tau\ast \Psi (t,\cdot)\|_{L^p}^p dt\frac {d\tau} \tau\notag\\
\leq &C\int_0^T \|\Psi (t,\cdot)\|_{\dot B_p^{\frac \alpha p,p}}^p dt=C\int_0^T \|(-\Delta)^{\frac 12-\frac \alpha2}\delta_{(\mathbb A_{[\theta]} \theta)}(t,\cdot)\|_{\dot B_p^{\frac \alpha p,p}}^p dt,
\end{align*}
using again the Fubini theorem and \eqref{CAR_THERM}-\eqref{DEF_PSI} for the last two inequality. Thus,
\begin{align}\label{PREAL_GRONWALL_BESOV}
\|\delta_\theta\|_{L_t^p(\dot B_{p,x}^{\frac \alpha p,p})}^p  \leq&C\int_0^T \|\delta_{(\mathbb A_{[\theta]} \theta)}(t,\cdot)\|_{\dot B_p^{\frac \alpha p+1-\alpha,p}}^p dt.
\end{align}
In order to apply a Gr\"onwall type inequality we are thus led to control the Besov norm of the product of the differences. To this end we will use the following Kato-Ponce type inequality in homogeneous Besov spaces. From Theorem 2.1 in \cite{Naibo} we get that for $p\in [1,+\infty[ $ and $p\le p_1,p_2,\tilde p_1,\tilde p_2 $ such that $\frac 1p=\frac 1{p_1}+\frac 1{p_2}=\frac 1{\tilde p_1}+\frac 1{\tilde p_2} $ and $s>0$,
\begin{equation}
\label{KP_BESOV_SPACE_HOMO}
\|fg\|_{\dot B_{p}^{s,q}}\le C (\|f\|_{\dot B_{p_1}^{s,q}}\|g\|_{L^{p_2}}+\|g\|_{\dot B_{\tilde p_1}^{s,q}}\|f\|_{L^{\tilde p_2}}).
\end{equation}
Write now:
$$\|\delta_{(\mathbb A_{[\theta]} \theta)}(t,\cdot)\|_{\dot B_p^{\frac \alpha p+1-\alpha,p}}=\|[{\mathbb A_{[\theta_1]}\theta_1-\mathbb A_{[\theta_2]} \theta_2}](t,\cdot)\|_{\dot B_p^{\frac \alpha p+1-\alpha,p}},$$
and by the linearity of the quantity $\mathbb A_{[\theta]}$ we have
$$\|\delta_{(\mathbb A_{[\theta]} \theta)}(t,\cdot)\|_{\dot B_p^{\frac \alpha p+1-\alpha,p}}\le \|(\mathbb A_{[\theta_1]}-\mathbb A_{[\theta_2]})\theta_1(t,\cdot)\|_{\dot B_p^{\frac \alpha p+1-\alpha,p}}+\|\mathbb A_{[\theta_2]}(\theta_1-\theta_2)(t,\cdot)\|_{\dot B_p^{\frac \alpha p+1-\alpha,p}}$$
thus, using the estimate (\ref{KP_BESOV_SPACE_HOMO}) above we obtain
\begin{eqnarray*}
&\le& C\Big(\|(\mathbb A_{[\theta_1]}-\mathbb A_{[\theta_2]})(t,\cdot)\|_{\dot B_{ p_1}^{\frac \alpha p+1-\alpha,p}}\|\theta_1\|_{L^{ p_2}}+\|(\mathbb A_{[\theta_1]}-\mathbb A_{[\theta_2]})(t,\cdot)\|_{L^{ p_2}}\|\theta_1\|_{\dot B_{ p_1}^{\frac \alpha p+1-\alpha,p}}\\
&&+\|(\theta_1-\theta_2)(t,\cdot)\|_{\dot B_{ p_1}^{\frac \alpha p+1-\alpha,p}}\|\mathbb A_{[\theta_2]}(t,\cdot)\|_{L^{ p_2}}+\|(\theta_1-\theta_2)(t,\cdot)\|_{L^{ p_2}}\|\mathbb A_{[\theta_2]}(t,\cdot)\|_{\dot B_{ p_1}^{\frac \alpha p+1-\alpha,p}}\Big).
\end{eqnarray*}
Recalling the stability of $\mathbb A $ in Besov and $L^p$ norms (see \eqref{HypoLPBorne1} and \eqref{HypoBesovBorne01}) we derive:
\begin{align*}
\|\delta_{(\mathbb A_{[\theta]} \theta)}(t,\cdot)\|_{\dot B_p^{\frac \alpha p+1-\alpha,p}}\le& C \Big(\|(\theta_1-\theta_2)(t,\cdot)\|_{\dot B_{ p_1}^{\frac \alpha p+1-\alpha,p}}(\|\theta_1(t,\cdot)\|_{L^{ p_2}}+\|\theta_2(t,\cdot)\|_{L^{ p_2}})\\
&+\|(\theta_1-\theta_2)(t,\cdot)\|_{L^{ p_2}}(\|\theta_1(t,\cdot)\|_{\dot B_{ p_1}^{\frac \alpha p+1-\alpha,p}}+\|\theta_2(t,\cdot)\|_{\dot B_{ p_1}^{\frac \alpha p+1-\alpha,p}})\Big).
\end{align*}
Plugging this bound into \eqref{PREAL_GRONWALL_BESOV} and using the maximum principle (\ref{MaximumEstimate}) yields:
 \begin{align}\label{PREAL_GRONWALL_BESOV_2}
\|\delta_\theta\|_{L_t^p(\dot B_{p,x}^{\frac \alpha p,p})}^p  \leq&C\Big(\int_0^T \|\delta_\theta(t,\cdot)\|_{\dot B_{p_1}^{\frac \alpha p+1-\alpha,p}}^p(2\|\theta_0\|_{L^{p_2}})^p dt\notag\\
&+\int_0^T \|\delta_\theta(t,\cdot)\|_{L^{p_2}}^p(\|\theta_1(t,\cdot)\|_{\dot B_{p_1}^{\frac \alpha p+1-\alpha,p}}^p+\|\theta_2(t,\cdot)\|_{\dot B_{p_1}^{\frac \alpha p+1-\alpha,p}}^p) dt\Big).
\end{align}
We also recall the following embedding between homogeneous Besov spaces (see Theorem A in \cite{han:95} or  (1.1) in  \cite{trie:14}   and Proposition 2.2 in \cite{sawa:18} in the inhomogeneous case see):
\begin{equation*}
\dot B^{s_0}_{p_0,q_0}\hookrightarrow \dot B^{s_1}_{p_1,q_1}\,\text{for}\,p_0,p_1,q_0,q_1 \in[1,+\infty],\,q_0\le q_1,\,p_0\le p_1, s_0>s_1 \,\text{such that}\,s_0-n/p_0= s_1-n/p_1. 
\end{equation*}
Apply this estimate with $ s_0=\frac\alpha p$, $s_1=\frac \alpha p+1-\alpha$ and $p_0=p$, which gives $s_0-s_1=\alpha-1=\varepsilon=n(\frac 1{p_0}-\frac 1{p_1})\iff \frac1{p_1}=\frac{1}{p_0}-\frac{\varepsilon}{n} $ and $q_0=q_1=p $, to derive in \eqref{PREAL_GRONWALL_BESOV_2}
 \begin{align}\label{PREAL_GRONWALL_BESOV_3}
\|\delta_\theta\|_{L_t^p(\dot B_{p,x}^{\frac \alpha p,p})}^p  \leq&C\int_0^T \|\delta_\theta(t,\cdot)\|_{\dot B_{p}^{\frac \alpha p,p}}^p(2\|\theta_0\|_{L^{p_2}})^p dt+C\int_0^T \|\delta_\theta(t,\cdot)\|_{L^{p_2}}^p(\|\theta_1(t,\cdot)\|_{\dot B_{p}^{\frac \alpha p,p}}^p+\|\theta_2(t,\cdot)\|_{\dot B_{p}^{\frac \alpha p,p}}^p) dt\notag\\
\leq &C\|\theta_0\|_{L^{p_2}}^p \int_0^T \|\delta_\theta(t,\cdot)\|_{\dot B_{p}^{\frac \alpha p,p}}^p dt +C\|\delta_\theta\|_{L_t^\infty(L_x^{p_2})}^p\int_0^T(\|\theta_1(t,\cdot)\|_{\dot B_{p}^{\frac \alpha p,p}}^p+\|\theta_2(t,\cdot)\|_{\dot B_{p}^{\frac \alpha p,p}}^p) dt\notag \\
\leq &C\|\theta_0\|_{L^{p_2}}^p \int_0^T \|\delta_\theta(t,\cdot)\|_{\dot B_{p}^{\frac \alpha p,p}}^p dt +C\|\delta_\theta\|_{L_t^\infty(L_x^{p_2})}^p(\|\theta_1\|_{L^p_t(\dot B_{p,x}^{\frac \alpha p,p})}^p+\|\theta_2\|_{L^p_t(\dot B_{p}^{\frac \alpha p,p})}^p)\notag \\ 
\le &C \|\delta_\theta\|_{L_t^\infty(L_x^{p_2})}^p,
\end{align}
where, for the last inequality, we use the Gr\"onwall lemma (recall the energy estimates (\ref{AprioriBesovIntro}), the maximum principle (\ref{MaximumEstimate}) and the information (\ref{ControlNormeDonneeInitiale})). Recall as well that we have  $\frac{1}{p}=\frac{1}{p_1}+\frac{1}{p_2}=\frac{1}{p}-\frac{\varepsilon}{n}+\frac{1}{p_2}\iff p_2=\frac{n}{\varepsilon} $ which is \textit{big} but finite in the subcritical case. Also, since we have assumed ${\bf \bar p}> \frac n\varepsilon $ the quantity $\|\delta_\theta\|_{L_t^\infty(L_x^{p_2})}^p $ is indeed finite (using again the maximum principle (\ref{MaximumEstimate}). To prove uniqueness it remains to justify that $\|\delta_\theta\|_{L_t^\infty(L_x^{p_2})}=0 $. Write for fixed $t$:
\begin{eqnarray}
\|\delta_\theta(t,\cdot)\|_{L^{p_2}}&\le& \left\|\int_{0}^{t} \nabla \mathfrak{p}^\alpha_{t-s}\ast  (\mathbb{A}_{[\theta_1]}\,\theta_1-\mathbb{A}_{[\theta_2]}\,\theta_2)(s,\cdot)ds\right\|_{L^{p_2}}\notag\\
&\le &\int_{0}^{t} \left\|\nabla \mathfrak{p}^\alpha_{t-s}\ast\delta_{(\mathbb{A}_{[\theta]}\,\theta)}(s,\cdot)\right\|_{L^{p_2}}ds\le \int_0^t \|\nabla \mathfrak{p}^\alpha_{t-s} \|_{L^{r_1}}\|\delta_{(\mathbb{A}_{[\theta]}\,\theta)}(s,\cdot)\|_{L^{r_2}}ds,\label{TO_UNIQUENESS_IN_L_INFTY_L_P2}
\end{eqnarray}
with $1+\frac{1}{p_2}= \frac{1}{r_1}+\frac{1}{r_2}$ with $r_1>1 $ so that we still have some margin to apply the H\"older inequality for the contribution $ \|\delta_{(\mathbb{A}_{[\theta]}\,\theta)}(s,\cdot)\|_{L^{r_2}}$ in which $r_2<p_2 $. From Lemma \ref{LEMME_HK} in the appendix, observe that 
$$\|\nabla \mathfrak{p}^\alpha_{t-s}\|_{L^{r_1}}\le C(t-s)^{-\frac 1\alpha-\frac{n}{\alpha r_1'}},\ \frac1{r_1}+\frac1{r_1'}=1,$$
which for $\alpha=1+\varepsilon $ gives an integrable singularity provided $r_1 $ is sufficiently close to 1. This is in fact equivalent to $\frac 1\alpha+\frac{n}{\alpha r_1'}<1$, from which we easily deduce the condition $\frac{n}{\varepsilon}<r'_1$ (since $\alpha=1+\varepsilon$). Now, with $\frac{1}{r_2}=\frac{1}{p_2}+\frac{1}{q_2}$
\begin{align*}
\|\delta_{(\mathbb{A}_{[\theta]}\,\theta)}(s,\cdot)\|_{L^{r_2}}&\le \|(\mathbb A_{[\theta_1]}-\mathbb A_{[\theta_2]}) \theta_1(s,\cdot)\|_{L^{r_2}}+\|(\theta_1-\theta_2)\mathbb A_{[\theta_2]}(s,\cdot)\|_{L^{r_2}}\\
&\le \|(\mathbb A_{[\theta_1]}-\mathbb A_{[\theta_2]})(s,\cdot)\|_{L^{p_2}} \|\theta_1(s,\cdot)\|_{L^{q_2}}+\|(\theta_1-\theta_2)(s,\cdot)\|_{L^{p_2}} \|\mathbb A_{[\theta_2]}(s,\cdot)\|_{L^{q_2}},
\end{align*}
with $q_2=\frac{r_2p_2}{p_2-r_2}$. Still from the stability of $\mathbb A $  (see \eqref{HypoLPBorne1} and \eqref{HypoBesovBorne01}) we get:
\begin{align*}
\|\delta_{(\mathbb{A}_{[\theta]}\,\theta)}(s,\cdot)\|_{L^{r_2}}
&\le \|(\theta_1-\theta_2)(s,\cdot)\|_{L^{p_2}} (\|\theta_1(s,\cdot)\|_{L^{q_2}}+\|\theta_2(s,\cdot)\|_{L^{q_2}})\\
&\le \|\delta_\theta\|_{L_t^\infty (L_x^{p_2})}(\|\theta_1\|_{L_t^\infty (L_x^{q_2})}+\|\theta_2\|_{L_t^\infty (L_x^{q_2})}),
\end{align*}
which plugged into \eqref{TO_UNIQUENESS_IN_L_INFTY_L_P2} yields for $t\in [0,T]$ with fixed $T>0 $,
\begin{align*}
 \|\delta_\theta(t,\cdot)\|_{L^{p_2}}\le CT^\zeta \|\delta_\theta\|_{L_t^\infty(L_x^{p_2})}(\|\theta_1\|_{L_t^\infty (L_x^{q_2})}+\|\theta_2\|_{L_t^\infty(L_x^{q_2})}), 
\end{align*}
for some exponent $\zeta >0$. By the maximum principle (\ref{MaximumEstimate}), the right-hand side above is bounded as long as we have $p_2, q_2\leq \bar{\mathbf p}$. Observe that the relationships $1+\frac{1}{p_2}= \frac{1}{r_1}+\frac{1}{r_2}$ and $\frac{1}{r_2}=\frac{1}{p_2}+\frac{1}{q_2}$ imply that $q_2=r'_1$, thus by the previous condition $\frac{n}{\varepsilon}<r'_1$ we need $\frac{n}{\varepsilon}<q_2$. Finally, since $p_2=\frac{n}{\varepsilon}$, we deduce the condition $ \frac{n}{\varepsilon}<\bar{\mathbf p}$. Thus, we obtain the uniqueness in $L_t^\infty (L_x^{p_2})$ if $T $ is small enough. From \eqref{PREAL_GRONWALL_BESOV_3} we eventually derive uniqueness in $L^p_t(\dot B_{p,x}^{\frac \alpha p, p})\underset{q\in [2,\bar {\mathbf p}]}{\bigcap}L_t^\infty (L_x^q)$. We again emphasize that for a fixed $\alpha=1+\varepsilon>0 $ a spatial integrability of order  $\bar {\mathbf p}:=\bar {\mathbf p}(\varepsilon) $ which goes to infinity when $\varepsilon $ goes to zero, is needed to guarantee uniqueness of the weak solutions in the considered function class. \hfill $\blacksquare $\\

\subsection{From energy estimates to uniform in time controls}
In order to prove our main Theorem \ref{THM_HOLDER} we need to perform a previous step which will prepare the information available on the solutions. Indeed, by estimate (\ref{ControlBesovApriori1}) of Theorem \ref{Theoreme1} we have for the weak solutions of equation (\ref{EquationPrincipaleIntro}) the information $\theta\in  L^{p}([0,T], \dot{B}^{\frac{\alpha}{p},p}_p(\mathbb{R}^n))$ with ${\bf 2}\leq p\leq{\bf\bar{p}}<+\infty$ and for some $0<T<+\infty$, but this fact is not accurate enough for our purposes as we need some bounded in time ($L^\infty_t$)-regularity control. To establish such a result, we will start from a mild representation of the weak solutions to \eqref{EquationPrincipaleIntro}. 
In this sense we have the following results.
\begin{Proposition}[From $L^\infty_t(L^p_x)$- to $L^\infty_t(L^1_x)$]\label{Propo_LinftyL1} 
Let the assumptions of Lemma \ref{LEMME_MILD} be in force and $\theta $ be  a weak solution to the nonlinear equation \eqref{EquationPrincipaleIntro}. Then this weak solution $\theta$ belongs to the space $L^{\infty}([0,T], L^1(\mathbb{R}^n))$ and thus by interpolation, it belongs to the spaces $L^{\infty}([0,T], L^p(\mathbb{R}^n))$ with $1<p<\bf \bar p$ and where $0<T<+\infty$.
\end{Proposition}
{\bf Proof.} Start from the mild representation \eqref{REP_MILD_W_SOL}. By the Young inequalities for the convolution, one obtains
\begin{eqnarray*}
\|\theta(t,\cdot)\|_{L^1}&\leq &\|\mathfrak{p}^\alpha_t\ast \theta_0\|_{L^1}+\left\| \int_{0}^{t}\mathfrak{p}^\alpha_{t-s}\ast\nabla\cdot(\mathbb{A}_{[\theta]}\,\theta(s,\cdot))ds\right\|_{L^1}\\
&\leq &\|\mathfrak{p}^\alpha_t\|_{L^1} \|\theta_0\|_{L^1}+ \int_{0}^{t}\|\nabla \mathfrak{p}^\alpha_{t-s}\|_{L^1}\|\mathbb{A}_{[\theta]}\,\theta(s,\cdot)\|_{L^1}ds\\
&\leq & \|\theta_0\|_{L^1}+ C\int_{0}^{t}(t-s)^{-\frac{1}{\alpha}}\|\mathbb{A}_{[\theta]}\|_{L^2}\|\theta(s,\cdot)\|_{L^2}ds\\
&\leq & \|\theta_0\|_{L^1}+ C\|\mathbb{A}_{[\theta]}\|_{L^\infty_t(L^2_x)}\|\theta\|_{L^\infty_t(L^2_x)}\int_{0}^{t}(t-s)^{-\frac{1}{\alpha}}ds,
\end{eqnarray*}
but since $\alpha=1+\varepsilon>1$, the integral above is bounded and since we have by (\ref{HypoLPBorne1}) and by the maximum principle (\ref{MaximumPrinciple1}) the estimates 
$$\|\mathbb{A}_{[\theta]}\|_{L^\infty_t(L^2_x)}\leq C\|\theta\|_{L^\infty_t(L^2_x)}\leq C\|\theta_0\|_{L^2}<+\infty,$$
we obtain the following control
\begin{equation}\label{EstimationL1NonLineaire}
\|\theta(t,\cdot)\|_{L^1}\leq \|\theta_0\|_{L^1}+C\|\theta_0\|_{L^2}^2t^{1-\frac{1}{\alpha}},
\end{equation}
which is bounded for $0\leq t\leq T<+\infty$ and we finally deduce that the weak solution $\theta$ belongs to the space $L^{\infty}([0,T], L^1(\mathbb{R}^n))$. \hfill $\blacksquare$\\
\begin{Remarque}
Note that if $1<p<2$, then by interpolation, for some $0<\nu<1$ such that $\frac1p=1-\frac{\nu}2$ we have
\begin{equation}\label{EstimationLpinterpolation}
\|\theta\|_{L^\infty_t(L^p_x)}\leq \|\theta\|_{L^\infty_t(L^1_x)}^{1-\nu}\|\theta\|_{L^\infty_t(L^2_x)}^\nu\leq C\big(\|\theta_0\|_{L^1}+\|\theta_0\|_{L^2}^2\big)^{1-\nu}\|\theta_0\|^\nu_{L^2},
\end{equation}
where in the last estimate we used (\ref{EstimationL1NonLineaire}) and the maximum principle (\ref{MaximumPrinciple1}). Note in particular that if $\mu\geq 1$, using the hypothesis \A{D} we have
\begin{equation}\label{EstimationLpinterpolation1}
\|\theta\|_{L^\infty_t(L^p_x)}\leq C\mu^2.
\end{equation}
Note that $C=C(T)$ but we will work with a fixed $0<T<+\infty$.
\end{Remarque}
In the next result we will show how to obtain a small gain of regularity for weak solutions. 
\begin{Proposition}\label{PropoLpLinfty1}
Let the dimension $n\geq 2$ and consider over $\mathbb{R}^{n}$ the equation (\ref{EquationPrincipaleIntro}) where the diffusion operator $\mathcal{L}^{\alpha}$ satisfies \A{A} with $\alpha=1+\varepsilon$, the nonlinear drift $\mathbb{A}_{[\theta]}$ defined in \A{B} enjoys the stability properties of \A{C} and the initial data $\theta_{0}$ satisfies \A{D}. Consider an associated weak solution $\theta(t,x)$ to the nonlinear equation \eqref{EquationPrincipaleIntro}. Then this weak solution $\theta$ belongs to the space $L^{\infty}([1,N], \dot{W}^{\frac{\sigma_0}{2},p_0}(\mathbb{R}^n))$ where $1<p_0<2$, with $N\geq 10$ and 
\begin{equation}\label{Lower_ConditionAlpha1}
\sigma_0=1+3\varepsilon= \alpha+2\varepsilon\leq 2\alpha.
\end{equation}
Moreover we have the estimate
\begin{equation}\label{EstimationLpBesovLinftSobo}
\|\theta\|_{L^{\infty}_t (\dot{W}^{\frac{\sigma_0}{2},p_0}_x)}\leq C\Big(\|\theta_0\|_{L^{{p_0}'}}\|\theta_0\|_{L^2}+\big(\|\theta_0\|_{L^1}+\|\theta_0\|_{L^2}^2\big)^{1-\nu}\|\theta_0\|^\nu_{L^2}\Big)<+\infty.
\end{equation}
\end{Proposition}
{\bf Proof.} We introduce a positive, smooth cut-off in time function $\phi \in \mathcal{C}^\infty_0(\mathbb{R})$ such that $\phi(s)\equiv0$ over $[0,\tfrac12]\cap [N+1, +\infty[$, such that $\phi(s)\equiv1$ over $[1, N]$ and such that $\|\phi\|_{L^\infty}=1$ and we define for $(t,x)\in \R_+\times \R^n $, the function 
$$u(t,x)=\phi(t)\theta(t,x).$$
Note that the functions $u(t,x)$ and $\theta(t,x)$ coincide if $1\leq t\leq N$ and moreover $u$ satisfies the following equation
$$
\begin{cases}
\partial_t u-\nabla\cdot(\mathbb{A}_{[\theta]}\,u)+\mathcal{L}^{\alpha}u-(\partial_t \phi) \theta=0,\qquad div(\mathbb{A}_{[\theta]})=0, \qquad 0<\alpha<2,\\[2mm]
u(t,x)=0,\quad t\in [0,\tfrac12].
\end{cases}
$$
Similarly to Lemma \ref{LEMME_MILD}, and recalling that  $(\mathfrak{p}^\alpha)_{t\ge 0} $ stands for the semi-group kernel associated with the operator $ {\mathcal L}^{\alpha}$, it holds from the Duhamel representation formula that for any $t_0\in [0,\tfrac12]$:
\begin{equation*}
u(t,x)= \int_{t_0}^{t}\mathfrak{p}^\alpha_{t-s}\ast\nabla\cdot(\mathbb{A}_{[\theta]}\,\theta(s,x)\phi(s))ds+\int_{t_0}^{t}\mathfrak{p}^\alpha_{t-s}\ast\theta(s,x)\partial_s\phi(s)ds.
\end{equation*}
Since $\alpha=1+\varepsilon$ for some small $0<\varepsilon<1$. The point is now to investigate, for a fixed $1<t<N$, the Sobolev $\dot W^{\frac{\sigma_0}{2}, p_0} $-norm in the space variable of the function $u(t, \cdot)$ where $\sigma_0$ satisfies (\ref{Lower_ConditionAlpha1}) and where the parameter $p_0$ satisfies $1<p_0<2$.\\

Indeed, write first for fixed $t$ and $t_0 $:
\begin{eqnarray*}
\|u(t,\cdot)\|_{\dot W^{\frac{\sigma_0}{2}, p_0}}&\le& \left\|\int_{t_0}^{t}(-\Delta)^{\frac{\sigma_0}{4}}\mathfrak{p}^\alpha_{t-s}\ast\nabla\cdot (\mathbb{A}_{[\theta]}\,\theta(s,\cdot)\phi(s))ds\right\|_{L^{p_0}}+\left\|\int_{t_0}^{t}(-\Delta)^{\frac{\sigma_0}{4}}\mathfrak{p}^\alpha_{t-s}\ast\theta(s,\cdot)\partial_s\phi(s)ds\right\|_{L^{p_0}}.
\end{eqnarray*}
Next we integrate with respect to $t_0\in [0,\frac{1}{2}]$ and we obtain:
\begin{eqnarray*}
\left(\int_0^{\frac12 }\|u(t,\cdot)\|_{\dot W^{\frac{\sigma_0}{2}, {p_0}}}^2 dt_0\right)^{\frac 12} &\le&\left\|\int_{\cdot}^{t}(-\Delta)^{\frac{1}{2}}\mathfrak{p}^\alpha_{t-s}\ast (-\Delta)^{\frac{\sigma_0} 4}(\mathbb{A}_{[\theta]}\,\theta(s,\cdot)\phi(s))ds\right\|_{L_{t_0}^2(L_x^{p_0})}\\
&&+\left\|\int_{\cdot}^{t}(-\Delta)^{\frac{\sigma_0}{4}}\mathfrak{p}^\alpha_{t-s}\ast\theta(s,\cdot)\partial_s\phi(s)ds\right\|_{L_{t_0}^2(L_x^{p_0})}.
\end{eqnarray*}
Since by (\ref{Lower_ConditionAlpha1}) we have $\sigma_0=1+3\varepsilon$ and thus we can write 
$$(-\Delta)^{\frac{\sigma_0}4}=(-\Delta)^{\frac{1+3\varepsilon} 4}=(-\Delta)^{\frac{1+\varepsilon} 4}(-\Delta)^{\frac{\varepsilon} 2}=(-\Delta)^{\frac{\alpha} 4}(-\Delta)^{\frac{\varepsilon} 2},$$
and it yields
\begin{eqnarray*}
&&\le\left\|\int_{\cdot}^{t}(-\Delta)^{\frac{1}{2}}\mathfrak{p}^\alpha_{t-s}\ast (-\Delta)^{\frac{\alpha} 4}(-\Delta)^{\frac{\varepsilon} 2}(\mathbb{A}_{[\theta]}\,\theta(s,\cdot)\phi(s))ds\right\|_{L_{t_0}^2(L_x^{p_0})}+\left\|\int_{\cdot}^{t}(-\Delta)^{\frac{\sigma_0}{4}}\mathfrak{p}^\alpha_{t-s}\ast\theta(s,\cdot)\partial_s\phi(s)ds\right\|_{L_{t_0}^2(L_x^{p_0})}\\
&&\le\left\|\int_{\cdot}^{t}(-\Delta)^{\frac{1+\varepsilon}{2}}\mathfrak{p}^\alpha_{t-s}\ast (-\Delta)^{\frac{\alpha} 4}(\mathbb{A}_{[\theta]}\,\theta(s,\cdot)\phi(s))ds\right\|_{L_{t_0}^2(L_x^{p_0})}+\left\|\int_{\cdot}^{t}(-\Delta)^{\frac{\sigma_0}{4}}\mathfrak{p}^\alpha_{t-s}\ast\theta(s,\cdot)\partial_s\phi(s)ds\right\|_{L_{t_0}^2(L_x^{p_0})}\\
&&\le\left\|\int_{\cdot}^{t}(-\Delta)^{\frac{\alpha}{2}}\mathfrak{p}^\alpha_{t-s}\ast (-\Delta)^{\frac{\alpha} 4}(\mathbb{A}_{[\theta]}\,\theta(s,\cdot)\phi(s))ds\right\|_{L_{t_0}^2(L_x^{p_0})}+\left\|\int_{\cdot}^{t}(-\Delta)^{\frac{\sigma_0}{4}}\mathfrak{p}^\alpha_{t-s}\ast\theta(s,\cdot)\partial_s\phi(s)ds\right\|_{L_{t_0}^2(L_x^{p_0})}.
\end{eqnarray*}
Now, we can apply  maximal-regularity type inequalities  to the kernels $(-\Delta)^{\frac{\alpha} 2}\mathfrak{p}^\alpha_{t-s}$ and $(-\Delta)^{\frac{\sigma_0}{4}}\mathfrak{p}^\alpha_{t-s}$ (since $\sigma_0\leq 2\alpha$ by (\ref{Lower_ConditionAlpha1})), to obtain
\begin{equation}\label{E1_E2}
\left(\int_0^{\frac12 }\|u(t,\cdot)\|_{\dot W^{\frac{\sigma_0}{2}, {p_0}}}^2 dt_0\right)^{\frac 12} \leq C\Big(\underbrace{\left\| (-\Delta)^{\frac{\alpha} 4}(\mathbb{A}_{[\theta]}\,\theta\phi)\right\|_{L_{t_0}^2(L_x^{p_0})}}_{E_1}+\underbrace{\left\|\theta\partial_s\phi\right\|_{L_{t_0}^2(L_x^{p_0})}}_{E_2}\Big).
\end{equation}
For the term $E_1$ above we apply the Leibniz fractional rule (also known as the Kato-Ponce inequality):
Namely, for $s>0 $, $\mathfrak{r}>1 $ and  $p_1,q_1,p_2,q_2 \in]1,+\infty]$ such that $\frac {1}{\mathfrak{r}}=\frac 1{p_1}+\frac1{q_1} =\frac 1{p_2}+\frac1{q_2} $ it holds that for two functions $f,g:\R^n\longrightarrow \R $,
\begin{equation}\label{LEIBNIZ}
\|(-\Delta)^{\frac s2}(fg)\|_{L^\mathfrak{r}}\le C\Big(\|f\|_{L^{p_1}}\|(-\Delta)^{\frac s 2}g\|_{L^{q_1}} +\|(-\Delta)^{\frac s 2}f\|_{L^{p_2}}\|g\|_{L^{q_2}} \Big).
\end{equation}
Inequality \eqref{LEIBNIZ} is a particular case of estimate (1) in \cite{GrafakosOh}. See also \cite{Naibo} for several versions of the above fractional Leibniz derivation rule.  Writing $\frac 1{p_0}=\frac 12+\frac 1{{p_0}'}$ (where $1<{p_0}<2$ can be very close to $2$ and $2<{p_0}'\leq {\bf \bar p}$ can be potentially very large) and taking $p_1={p_0}', q_1=2, p_2=2, q_2={p_0}'$, we then derive from 
\eqref{LEIBNIZ} and from the definition of $E_1$ given in  \eqref{E1_E2} that:
\begin{eqnarray*}
E_1&=&\left\| (-\Delta)^{\frac{\alpha} 4}(\mathbb{A}_{[\theta]}\,\theta\phi)\right\|_{L_{t_0}^2(L_x^{p_0})}\\
&\leq &C\Big(\left\| \|(-\Delta)^{\frac{\alpha} 4}\mathbb{A}_{[\theta]}\|_{L^2_x}\|\theta\|_{L^{{p_0}'}_x}\phi\right\|_{L_{t_0}^2}+\left\| \|\mathbb{A}_{[\theta]}\|_{L^{{p_0}'}_x}\|(-\Delta)^{\frac{\alpha} 4}\theta\|_{L^{2}_x}\phi\right\|_{L_{t_0}^2}\Big),
\end{eqnarray*}
since $\phi\in L^\infty$ and $\|\theta\|_{L^\infty_t(L^{{p_0}'}_x)}<+\infty$ we have
\begin{eqnarray*}
E_1&\leq &C\|\phi\|_{L^\infty}\Big(\left\| \|(-\Delta)^{\frac{\alpha} 4}\mathbb{A}_{[\theta]}\|_{L^2_x}\right\|_{L_{t_0}^2}\|\theta\|_{L^\infty_tL^{{p_0}'}_x}+\left\| \|(-\Delta)^{\frac{\alpha} 4}\theta\|_{L^{2}_x}\right\|_{L_{t_0}^2}\|\mathbb{A}_{[\theta]}\|_{L^\infty_tL^{{p_0}'}_x}\Big).
\end{eqnarray*}
Recalling now the estimates $\|\phi\|_{L^\infty}=1$ as well as the controls (as a consequence of (\ref{HypoBesovBorne02})-(\ref{HypoLPBorne1}) and (\ref{MaximumPrinciple1}))
$$ \|(-\Delta)^{\frac{\alpha} 4}\mathbb{A}_{[\theta]}\|_{L^2_x}\leq C \|(-\Delta)^{\frac{\alpha} 4}\theta\|_{L^2_x}\qquad \mbox{and }\quad \|\mathbb{A}_{[\theta]}\|_{L^\infty_t(L^{{p_0}'}_x)}\leq C\|\theta\|_{L^\infty_t(L^{{p_0}'}_x)},$$
we obtain the inequality 
$$E_1\leq C\|\theta\|_{L^\infty_t(L^{{p_0}'}_x)}\left\| \|(-\Delta)^{\frac{\alpha} 4}\theta\|_{L^{2}_x}\right\|_{L_{t_0}^2}=C\|\theta\|_{L^\infty_t(L^{{p_0}'}_x)}\|\theta\|_{L^2_t(\dot{H}^{\frac{\alpha}{2}}_x)}.$$
At this point we use the maximum principle (\ref{MaximumPrinciple1}) since $p_0'>2$:
$$\|\theta\|_{L^\infty_t(L^{{p_0}'}_x)}\leq C\|\theta_0\|_{L^{{p_0}'}},$$
and the a priori inequality (\ref{ControlBesovApriori1}) with $p=2$ (since we have the identification $\dot{B}^{\frac{\alpha}{2},2}_2(\R^n)=\dot{H}^{\frac{\alpha}{2}}(\R^n)$)
$$\|\theta\|_{L^2_t(\dot{H}^{\frac{\alpha}{2}}_x)}\leq \|\theta_0\|_{L^2},$$
to deduce that we have 
$$E_1\leq C\|\theta_0\|_{L^{{p_0}'}}\|\theta_0\|_{L^2}.$$
The term $E_2$ of (\ref{E1_E2}) is easier to study and we have, since $1<{p_0}<2$, by the expression (\ref{EstimationLpinterpolation}) above:
$$E_2\le C \|\theta\|_{L^\infty_t(L^{p_0}_x)}\|\partial_s\phi\|_{L^2_t}\leq C\big(\|\theta_0\|_{L^1}+\|\theta_0\|_{L^2}^2\big)^{1-\nu}\|\theta_0\|^\nu_{L^2}.$$
Thus, coming back to (\ref{E1_E2}) we deduce that 
$$\|u(t,\cdot)\|_{\dot W^{\frac{\sigma_0}{2}, {p_0}}}\leq C\|\theta_0\|_{L^{{p_0}'}}\|\theta_0\|_{L^2}+C\big(\|\theta_0\|_{L^1}+\|\theta_0\|_{L^2}^2\big)^{1-\nu}\|\theta_0\|^\nu_{L^2}.$$
but since $u$ and $\theta$ coincide over the interval $[1, N]$, we have that $\theta\in L^{\infty}([1,N], \dot W^{\frac{\sigma_0}{2},{p_0}})$ and 
\begin{equation}\label{ConclusionStep1}
\|\theta(t,\cdot)\|_{\dot W^{\frac{\sigma_0}{2},{p_0}}} \leq C\Big(\|\theta_0\|_{L^{{p_0}'}}\|\theta_0\|_{L^2}+\big(\|\theta_0\|_{L^1}+\|\theta_0\|_{L^2}^2\big)^{1-\nu}\|\theta_0\|^\nu_{L^2}\Big)<+\infty,
\end{equation}
and Proposition \ref{PropoLpLinfty1} is proven. \hfill $\blacksquare$
\begin{Remarque} Note that the gain of regularity $\sigma_0$ given in (\ref{Lower_ConditionAlpha1}) is independent from the value of $1<p_0<2$. 
Moreover, if $p_0\to 2$ then $p_0'\to +\infty$ and since we have $p_0'\leq {\bf \bar p}$ we should have ${\bf \bar p}\to +\infty$.
\end{Remarque}
\begin{Remarque}
Recall that by hypothesis \A{D} we have $\|\theta_0\|_{L^p}\leq \mu<+\infty$ for all $1\leq p\leq {\bf \bar p}$. Then the conclusion (\ref{ConclusionStep1}) of the previous proposition can be rewritten as 
$$\|\theta(t,\cdot)\|_{\dot W^{\frac{\sigma_0}{2},{p_0}}} \leq C\Big(\mu^2+\big(\mu+\mu^2\big)^{1-\nu}\mu^\nu\Big),$$
for $1\leq t\leq N$ and since $\mu\geq 1$ we have
\begin{equation}\label{ConclusionStep2}
\|\theta(t,\cdot)\|_{\dot W^{\frac{\sigma_0}{2},{p_0}}} \leq C\mu^2.
\end{equation}
\end{Remarque}
\begin{Corollaire}[From $L^p$-Besov to $L^\infty$-Besov]\label{CorollaireLinftyBesov}
\begin{itemize}
\item[]
\item[1)]For some $0<\varepsilon_0<\varepsilon$ such that $\frac{n(2-p_0)}{2p_0}=\varepsilon-\varepsilon_0$, we have $\theta\in L^{\infty}([1,N], \dot H^{\frac \alpha 2+\varepsilon_0})$, in particular we have
$$\|\theta(t,\cdot)\|_{\dot H^{\frac \alpha 2+\varepsilon_0}} \leq  C\mu^2.$$
\item[2)] For $0<\nu<1$, we set  $2<p_\nu<\bf\bar p$ by the condition $\frac{1}{p_\nu}=\frac{1-\nu}{ \bf\bar p}+\frac{\nu}{2}$. Then, for the same index $0<\varepsilon_0<\varepsilon$ as in the first point above, we have $\theta\in L^{\infty}([1,N], \dot B^{\nu(\frac{\alpha}{2}+\varepsilon_0),p_\nu}_{p_\nu})$, \emph{i.e.} for $1\leq t\leq N$:
$$\|\theta(t,\cdot)\|_{\dot B^{\nu(\frac{\alpha}{2}+\varepsilon_0),p_\nu}_{p_\nu}}\leq C\mu^2.$$
\item[3)] For $1<p_0<2$ we have $\theta\in L^{\infty}([1,N], \dot B^{\frac{\alpha}{2}+\varepsilon,p_0}_{2})$:
$$\|\theta(t,\cdot)\|_{\dot B^{\frac{\alpha}{2}+\varepsilon,p_0}_{2}}\leq C\mu^2.$$
\end{itemize}
\end{Corollaire}
{\bf Proof.} For the first point, by the Sobolev embeddings if $\frac{\alpha}{2}+\varepsilon_0-\frac{n}{2}=\frac{\sigma_0}{2}-\frac{n}{p_0}$ and since we have $\sigma_0=\alpha+2\varepsilon$, we have for $1\leq t\leq N$:
\begin{equation}\label{EstimationSobolevMU2}
\|\theta(t,\cdot)\|_{\dot{H}^{\frac{\alpha}{2}+\varepsilon_0}}\leq C\|\theta(t,\cdot)\|_{\dot W^{\frac{\sigma_0}{2},p_0}}\leq C \mu^2.
\end{equation}
For the second point, by the complex interpolation theory we have, for $0<\nu<1$ and for $1\leq t\leq N$:
$$\left[L^\infty_t(L^{\bf \bar  p}_x)\; ;\; L^{\infty}_t(\dot H^{\frac{\alpha}{2}+\varepsilon_0}_x) \right]_{\nu}= L^\infty_t(\dot W^{\nu(\frac{\alpha}{2}+\varepsilon_0),p_\nu}_x)\subset L^\infty_t(\dot B^{\nu(\frac{\alpha}{2}+\varepsilon_0),p_\nu}_{p_\nu,x}),$$
where $\frac{1}{p_\nu}=\frac{1-\nu}{ \bf\bar p}+\frac{\nu}{p_1}$ and since we have the space inclusion $\dot{W}^{s,p}\subset \dot{B}^{s,p}_{\max\{2,p\}}$ we can write
\begin{eqnarray*}
\|\theta\|_{L^\infty_t(\dot B^{\nu(\frac{\alpha}{2}+\varepsilon_0),p_\nu}_{p_\nu, x})}&\leq &\|\theta\|_{L^\infty_t(\dot W^{\nu(\frac{\alpha}{2}+\varepsilon_0),p_\nu}_x)}\leq C\|\theta\|_{L^\infty_t(L^{\bf \bar  p}_x)}^{1-\nu}\|\theta\|_{L^{\infty}_t(\dot H^{\frac{\alpha}{2}+\varepsilon_0}_x)}^\nu\\
&\leq &C(\mu)^{1-\nu}(\mu^2)^{\nu},
\end{eqnarray*}
where we used the hypothesis \A{D} for the first norm, the estimate (\ref{EstimationSobolevMU2}) for the second one and we finally obtain
$$\|\theta\|_{L^\infty_t(\dot B^{\nu(\frac{\alpha}{2}+\varepsilon_0),p_\nu}_{p_\nu, x})}\leq C\mu^2,$$
since we have $\mu\geq1$. For the last point, it is enough to recall that we have the space inclusion $\dot{W}^{\frac{\sigma_0}{2},p_0}\subset \dot{B}^{\frac{\sigma_0}{2},p_0}_{2}$ since $1<p_0<2$. Thus, by (\ref{ConclusionStep2}) we have for $1\leq t\leq N$
$$\|\theta(t,\cdot)\|_{\dot B^{\frac{\sigma_0}{2},p_0}_{2}}\leq \|\theta(t,\cdot)\|_{\dot W^{\frac{\sigma_0}{2},p_0}}\leq C\mu^2,$$
to conclude recall that by (\ref{Lower_ConditionAlpha1}) we have $\sigma_0=\alpha+2\varepsilon$.
\hfill $\blacksquare$

\begin{Remarque}\label{Remarque_EstimationBesovDrift}
Due to the hypotheses over the drift $\mathbb{A}$ stated in (\ref{HypoBesovBorne01}), we have for $1\leq t\leq N$:
\begin{itemize}
\item[1)] For $0<\nu\leq 1$ and for $2\leq p_\nu<\bf\bar p$ given by the condition $\frac{1}{p_\nu}=\frac{1-\nu}{ \bf\bar p}+\frac{\nu}{2}$:
\begin{equation}\label{EstimationLinftyBesovPgrandA}
\|\mathbb{A}_{[\theta]}(t,\cdot)\|_{\dot B^{\nu(\frac{\alpha}{2}+\varepsilon_0),p_\nu}_{p_\nu}}\leq C\|\theta(t,\cdot)\|_{\dot B^{\nu(\frac{\alpha}{2}+\varepsilon_0),p_\nu}_{p_\nu}}\leq C\mu^2.
\end{equation}
\item[2)] For $1<p_0<2$:
\begin{equation}\label{EstimationLinftyBesovPpetitA}
\|\mathbb{A}_{[\theta]}(t,\cdot)\|_{\dot B^{\frac{\alpha}{2}+\varepsilon,p_0}_{2}}\leq C\|\theta(t,\cdot)\|_{\dot B^{\frac{\alpha}{2}+\varepsilon,p_0}_{2}}\leq C\mu^2.
\end{equation}
\end{itemize}
\end{Remarque}
\mysection{H\"older Regularity}\label{SeccHolderRegularity}
To study the H\"older regularity of a solution $\theta(t,x)$ of equation (\ref{EquationPrincipaleIntro}), we will use the fact that the dual space of Hardy spaces $h^{\mathfrak{s}}(\mathbb{R}^n)$ are precisely\footnote{See \cite{Gold} for a proof of this fact and see \cite{Coifmann} and \cite{Stein2} for a detailed treatment on Hardy spaces.} the H\"older spaces $\mathcal{C}^{\gamma}(\mathbb{R}^n)$. Indeed, let $\frac{n}{n+1}<\mathfrak{s}<1$ and fix $\gamma$ by the relationship
\begin{equation}\label{FormuleGammaRegulariteHolderHardy}
0<\gamma=n(\frac{1}{\mathfrak{s}}-1)<1,
\end{equation}
then the dual of the local Hardy space $h^{\mathfrak{s}}(\mathbb{R}^n)$ is the Hölder space $\mathcal{C}^\gamma(\mathbb{R}^n)$, \emph{i.e.} we have the identification $(h^{\mathfrak{s}})'\simeq \mathcal{C}^\gamma$. 
We recall that the local Hardy space $h^{\mathfrak{s}}(\mathbb{R}^n)$ is the set of distributions $f$ that admits a \emph{molecular} decomposition of the form $f=\displaystyle{\sum_{j\in \mathbb{N}}}\lambda_j \psi_j$, where $(\lambda_j)_{j\in \mathbb{N}}$ is a sequence of numbers such that $\displaystyle{\sum_{j\in \mathbb{N}}}|\lambda_j|^\mathfrak{s}<+\infty$ and $(\psi_j)_{j\in \mathbb{N}}$ is a family of molecules given by the following definition.
\begin{Definition}\label{DefMolecules} Set $\frac{n}{n+1}<\mathfrak{s}<1$, define $\gamma$ by condition (\ref{FormuleGammaRegulariteHolderHardy}) and fix a real number $\omega$ such that $0<\gamma<\omega<1$. Consider a real parameter $\zeta>0$. For $0<r\ll 1$, we will say that an integrable function $\psi$ is a \emph{small molecule} with center $x_0\in \mathbb{R}^n$ and \emph{size} $\zeta r$ if we have
\begin{eqnarray}
& & \int_{\mathbb{R}^n} |\psi(x)||x-x_0|^{\omega}dx \leq   (\zeta r)^{\omega-\gamma}\mbox{, for } x_0\in \mathbb{R}^n \qquad \mbox{and }\qquad \|\psi\|_{L^\infty}  \leq  \frac{1}{ (\zeta r)^{n+\gamma}}\label{Hipo2},\\[1mm]
& &\int_{\mathbb{R}^n} \psi(x)dx=0\label{Hipo3}.
\end{eqnarray}
 In the case when $1\leq \zeta  r <+\infty$ (i.e. for big molecules), we only require conditions (\ref{Hipo2}) for the  molecule $\psi$ while the moment condition (\ref{Hipo3}) is dropped.
\end{Definition}
\begin{Remarque}\label{Remark200}
With the identification $(h^{\mathfrak{s}})'\simeq \mathcal{C}^\gamma$ and using the molecular decomposition of Hardy spaces, the fact that $\theta\in \mathcal{C}^\gamma$ is equivalent to the fact that $|\langle \theta, \psi\rangle_{\mathcal{C}^{\gamma}\times h^{\mathfrak{s}}}|<+\infty$ for all molecule $\psi$.
\end{Remarque}
\begin{Remarque}\label{Remark20}
The technical parameter $\omega$ satisfies the inequalities $0<\gamma<\omega$ and thus gives a maximum threshold for the H\"older regularity $\gamma$. We will always assume that we have 
\begin{equation}\label{RelationGammaOmegaAlpha1}
0<\gamma<\omega<1<\alpha=1+\varepsilon,
\end{equation}
where $\alpha$ is the smoothness degree associated to the diffusion operator $\mathcal{L}^{\alpha}$. These inequalities reflect the fact that it is not possible to obtain (by the method displayed here) a  H\"older regularity index $\gamma$ higher than the smoothness degree $\alpha$.
\end{Remarque}

\begin{Remarque}\label{Remark2}
Conditions (\ref{Hipo2}) imply the estimate $\|\psi\|_{L^1}\leq C\,  (\zeta r)^{-\gamma}$ (see e.g. Section 4.1, \textbf{3)} in \cite{DCHSM}). Thus, every molecule belongs to $L^q(\mathbb{R}^n)$ with $1<q<+\infty$ since we have
\begin{equation}\label{Hipo4}
\|\psi\|_{L^q} \leq C (\zeta  r)^{-n+\frac{n}{q}-\gamma}.
\end{equation}
Recall also that the Schwartz class $\mathcal{S}(\mathbb{R}^{n})$ is dense in $h^{\mathfrak{s}}(\mathbb{R}^{n})$, this fact is of course very useful for approximation procedures.
\end{Remarque}
We refer to \cite{DCHSM} for additional remarks about the former definition. See also \cite[Chapter III, \S 5.7]{Stein2}, \cite[Chapter XIV, \S 6.6]{Torchinski} or the article \cite{KN}.
\subsection{The dual equation}\label{Secc_EquationDuale}
Once we have described the elements of the Hardy spaces, we need to derive a dual equation from the original problem (\ref{EquationPrincipaleIntro}), and for this we proceed as follows. For $\mathcal{L}^{\alpha}$ with $\alpha=1+\varepsilon$ a L\'evy-type operator satisfying \A{A}, for a nonlinear drift $\mathbb{A}_{[\theta]}$ given by \A{B} and satisfying hypothesis \A{C} and for an initial data $\theta_{0}$ satisfying \A{D}, it follows from Theorem \ref{Theoreme1} that we can construct on the interval $[0,T]$, with $0<T<+\infty$ fixed, a corresponding weak solution $\theta(\cdot,\cdot)\in L^{\infty}\big([0,T], L^{p}(\mathbb{R}^{n})\big)\cap L^{p}([0,T], \dot{B}^{\frac{\alpha}{p},p}_p(\mathbb{R}^n))$ to \eqref{EquationPrincipaleIntro} that satisfies the inequalities (\ref{MaximumPrinciple1}) and (\ref{ControlBesovApriori1}).\\

Fix now a time $2< t\leq N$ with $N$ large enough (say $N\geq 10$, recall Proposition \ref{PropoLpLinfty1}), consider $\psi_0$ a molecule in the sense of the Definition \ref{DefMolecules} and consider a dual time variable $0\leq s\leq t$. The choice $t>2$ is here arbitrary and mainly performed for simplicity. We could have considered $t>t_0>0$ provided $ \mu:=\mu(t_0)$ is \textit{small} enough. With all these objects, to each molecule $\psi_0$ we can associate the following \emph{linear dual equation}
\begin{equation}\label{Evolution01}
\begin{cases}
\partial_s \psi(s,x) +\nabla\cdot [\mathbb{A}_{[\theta]}(t-s,x)\psi(s,x)]+\mathcal{L}^{\alpha}\psi(s,x)=0,\qquad s\in [0,t],\\[3mm]
\psi(0,x)=\psi_0(x),
\end{cases}
\end{equation}
where $\mathcal{L}^{\alpha}$ is a L\'evy-type operator of degree $\alpha=1+\varepsilon$ that satisfies hypothesis \A{A} and the drift $\mathbb{A}$ satifies the hypothesis \A{B} and \A{C}.\\

As said in the introduction, this equation shares many common features with the nonlinear equation (\ref{EquationPrincipaleIntro}). However, since we are in a linear setting, there are some particularities that have to be taken into account.  We state below the main results that are needed for the sequel (the proofs are given \cite{DCHSM} and \cite{DCHSM_POTA}). From the Proposition \ref{PropoLineaireExistence} to the Proposition \ref{Propo_Transfert} we state these results using for the sake of simplicity an initial data $\psi_0\in L^{1}(\R^{n})\cap L^{\infty}(\mathbb{R}^n)$.
\begin{Proposition}[Existence]\label{PropoLineaireExistence}
Let $n\geq 2$. If $\psi_0\in L^{1}(\R^{n})\cap L^{\infty}(\mathbb{R}^n)$ is a initial data, $\mathcal{L}^{\alpha}$ is a L\'evy-type operator that satisfy \A{A} with $0<\alpha<2$ and if the drift $\mathbb{A}_{[\theta]}$ satisfies the uniform boundedness conditions given in Remark \ref{Remarque_EstimationBesovDrift} (equations \eqref{EstimationLinftyBesovPgrandA} and \eqref{EstimationLinftyBesovPpetitA}) as well as the controls given in the hypothesis \A{C}, then there exists a weak solution $\psi(s,x)$ to equation \eqref{Evolution01} in $L^{\infty}([0,t], L^{q}(\mathbb{R}^n))$ with $1\leq q\leq +\infty$.\end{Proposition}
\begin{Proposition}[Maximum principle and Besov information]\label{PropoLineaireMaxetBesov}
Under the framework of Proposition \ref{PropoLineaireExistence}, weak solutions of equation (\ref{Evolution01}) satisify the following maximum principle for $s\in[0,t]$:
\begin{equation}\label{ControlLqAprioriMolecules1}
\|\psi(s,\cdot)\|_{L^q}\leq \|\psi_0\|_{L^q}\qquad \mbox{with}\quad 1\leq q\leq +\infty,
\end{equation}
and moreover we also have the following Besov a priori control for $2\leq q<+\infty$:
\begin{equation}\label{ControlBesovAprioriMolecules1}
\|\psi\|_{L^{q}_t(\dot{B}^{\frac{\alpha}{q},q}_{q,x})}\leq C \|\psi_0\|_{L^q}.
\end{equation}
\end{Proposition}
\begin{Proposition}[Positivity principle]\label{PropoLineairePositivite}
Let $\psi_{0}\in L^{1}(\R^{n})\cap L^{\infty}(\R^{n})$ be an initial data such that $0 \leq \psi_{0} \leq M$ a.e. where $M > 0$ is a constant. Under the hypotheses \A{A} for the operator $\mathcal{L}^{\alpha}$ and \A{B}-\A{C} for the drift $\mathbb{A}_{[\theta]}$, then the weak solution of equation (\ref{Evolution01}) satisfies $0\leq  \psi(s,x) \leq M$ for all $s\in [0,t]$.
\end{Proposition}
We can now state the following result (borrowed from \cite{KN}) which is crucial to prove Theorem \ref{THM_HOLDER} by using the Hardy-H\"older duality.
\begin{Proposition}[Transfer Property]\label{Propo_Transfert} Let $\theta_{0}$ be an initial data satisfying \A{D} and let $\theta(t,x)$ be a weak solution of the equation (\ref{EquationPrincipaleIntro}) in the interval $[0,T]$ where the L\'evy-type operator $\mathcal{L}^{\alpha}$ and the nonlinear drift $\mathbb{A}_{[\theta]}$ satisfy the hypotheses \A{A}, \A{B} and \A{C}. Let $2< t\leq N$ and let $\psi(s, x)$ be a solution with a molecular initial data $\psi_{0}$ of the backward problem (\ref{Evolution01}) for $0\leq s\leq t$. Then we have the identity
\begin{equation}\label{Formule_Preservation}
\int_{\mathbb{R}^n}\theta(t,x)\psi(0,x)dx=\int_{\mathbb{R}^n}\theta\left(t-s,x\right)\psi\left(s,x\right)dx.
\end{equation}
\end{Proposition}
\textbf{Proof.} We first consider the expression 
$$\partial_s\int_{\mathbb{R}^n}\theta(t-s,x)\psi(s,x)dx=\int_{\mathbb{R}^n}-\partial_t\theta(t-s,x)\psi(s,x)+\partial_s\psi(s,x)\theta(t-s,x)dx.$$
Using equations (\ref{EquationPrincipaleIntro}) and (\ref{Evolution01}) we obtain
{\small
\begin{eqnarray*}
\partial_s\int_{\mathbb{R}^n}\theta(t-s,x)\psi(s,x)dx&=&\int_{\mathbb{R}^n}\bigg[\textcolor{black}{-}\nabla\cdot\left(\mathbb{A}_{[\theta]}(t-s,x)\theta(t-s,x)\right)+\mathcal{L}^{\alpha}\theta(t-s,x)\bigg]\psi(s,x) \\[5mm]
&&+\bigg[ -\nabla\cdot\left((\mathbb{A}_{[\theta]}(t-s,x) \psi(s,x)\right)-\mathcal{L}^{\alpha}\psi(s,x) \bigg]\theta(t-s,x)dx
\end{eqnarray*}
\begin{eqnarray}
&=& \int_{\mathbb{R}^n} \bigg[\textcolor{black}{-}\nabla\cdot\left(\mathbb{A}_{[\theta]}(t-s,x)\theta(t-s,x)\right)\bigg]\psi(s,x)\textcolor{black}{+}\bigg[-\nabla\cdot\left(\mathbb{A}_{[\theta]}(t-s,x) \psi(s,x)\right)\bigg]\theta(t-s,x)dx\notag\\
&&+ \int_{\mathbb{R}^n}\left(\mathcal{L}^{\alpha}\theta(t-s,x)\right)\psi(s,x)- \left(\mathcal{L}^{\alpha}\psi(s,x) \right)\theta(t-s,x)dx.\label{Preservation}
\end{eqnarray}
}
From the symmetry of the operator $\mathcal{L}^{\alpha}$ we have
$$\int_{\mathbb{R}^n}\left(\mathcal{L}^{\alpha}\theta(t-s,x)\right)\psi(s,x)dx=\int_{\mathbb{R}^n} \left(\mathcal{L}^{\alpha}\psi(s,x) \right)\theta(t-s,x)dx,$$
and thus the second integral of the previous formula is null. Now, since the transport terms $\mathbb{A}_{[\theta]}$ is divergence free we obtain 
\begin{eqnarray*}
&&\int_{\mathbb{R}^n} \bigg[\nabla\cdot\left(\mathbb{A}_{[\theta]}(t-s,x)\theta(t-s,x)\right)\bigg]\psi(s,x)-\bigg[-\nabla\cdot\left(\mathbb{A}_{[\theta]}(t-s,x) \psi(s,x)\right)\bigg]\theta(t-s,x)dx\\
&=&\textcolor{black}{-}\int_{\mathbb{R}^n} \bigg[\mathbb{A}_{[\theta]}(t-s,x) -\mathbb{A}_{[\theta]}(t-s,x)\bigg]\theta(t-s,x)\textcolor{black}{\cdot} \nabla\psi(s,x)dx=0
\end{eqnarray*}
and this integral is null and all the expression (\ref{Preservation}) is equal to zero, so the integral quantity $\displaystyle{\int_{\mathbb{R}^n}\theta(t-s,x)\psi(s,x)dx}$, remains constant in time.  \hfill$\blacksquare$
\mysection{Proof of Theorem \ref{THM_HOLDER}.}\label{Secc_Theo2}
We want to establish the H\"older regularity $\mathcal{C}^{\gamma}(\R^n)$ for the solutions $\theta$ of the nonlinear equation (\ref{EquationPrincipaleIntro}) and following Remark \ref{Remark200} we need to prove that the duality bracket $|\langle \theta(t,\cdot),\psi_{0}(\cdot) \rangle_{\mathcal{C}^{\lambda}\times h^{\mathfrak{s}}}|$ is finite for all molecule $\psi_{0}$. Then by the transfer property given in Proposition \ref{Propo_Transfert} we have the identity 
\begin{equation*}
\langle \theta(t,\cdot),\psi_{0}(\cdot) \rangle_{\mathcal{C}^{\lambda}\times h^{\mathfrak{s}}}= \langle \theta(t-s,\cdot),\psi(s,\cdot) \rangle_{L^{p}\times L^{p'}},
\end{equation*}
which transforms a H\"older-Hardy bracket into a pure Lebesgue one for some ${\bf 2}\leq p\leq{\bf \bar p}$ and $\frac{1}{p}+\frac{1}{p'}=1$.
Applying the maximum principle (\ref{MaximumPrinciple1}) and using the hypothesis \A{D} in the previous identity we can deduce the following inequalities
\begin{eqnarray}
|\langle \theta(t,\cdot),\psi_{0}(\cdot) \rangle_{\mathcal{C}^{\lambda}\times h^{\mathfrak{s}}}|&\leq &\| \theta(t-s,\cdot)\|_{L^{p}}\| \psi(s,\cdot)\|_{L^{p'}}\nonumber\\
&\leq & \| \theta_0\|_{L^{p}}\| \psi(s,\cdot)\|_{L^{p'}}\leq \mu \|\psi(s,\cdot)\|_{L^{p'}},\label{Dualite}
\end{eqnarray}
where $\mu$ is given in (\ref{ControlNormeDonneeInitiale}). Thus, in order to prove Theorem \ref{THM_HOLDER}, we only need to estimate the quantity $\|\psi(s,\cdot)\|_{L^{p'}}$ that comes from a molecular initial data $\psi_{0}$.\\

Now, due to the maximum principle (applied this time to $\psi(s,x)$) we can divide our proof into two steps following the molecule's size. Indeed, for \emph{big} molecules, \emph{i.e.} if $\zeta r\geq C$, we have by (\ref{Hipo4}) the inequality
\begin{equation}\label{EstimationLpprimeMolecules0}
\|\psi(s,\cdot)\|_{L^{p'}}\leq\|\psi_0\|_{L^{p'}}\leq C (\zeta r)^{-n+\frac{n}{p'}-\gamma}<+\infty,
\end{equation}
which is immediately bounded. It only remains to study the $L^{p'}$ control for \textit{small} molecules and this is done in the following theorem:
\begin{Theoreme}\label{TheoLp'control}
Let the assumptions of Theorem \ref{THM_HOLDER} hold, $\psi_{0}$ be a small molecule and consider $\psi(s,\cdot)$ the associated solution of the backward problem (\ref{Evolution01}). There exists a small time $0<\mathcal{T}_{0}<1$ such that we have the following control of the $L^{p'}$-norm of $\psi(s,\cdot)$:
\begin{equation}\label{EstimationLpprimeMolecules1}
\|\psi(s,\cdot)\|_{L^{p'}}\leq C \mathcal{T}_{0}^{-n+\frac{n}{p'}-\gamma},\mathcal{T}_{0}\leq s\leq t-2,
\end{equation}
where $0<\gamma<\alpha$ and where $C>0$ is a positive constant.
\end{Theoreme}
{\bf Proof of Theorem \ref{THM_HOLDER}.}
Accepting for a while the previous Theorem \ref{TheoLp'control}, we have then a good control over the quantity $\|\psi(s,\cdot)\|_{L^{p'}}$ for \emph{big} and \emph{small} molecules and getting back to (\ref{Dualite}) we obtain that the duality bracket $|\langle \theta(t,\cdot), \psi_0\rangle_{\mathcal{C}^{\gamma}\times h^{\mathfrak{s}}}|$ is always bounded for any molecule $\psi_0$. This proves Theorem \ref{THM_HOLDER} by duality and thus the solutions $\theta(t,x)$ of the nonlinear equation (\ref{EquationPrincipaleIntro}) are $\gamma$-H\"older regular. \hfill $\blacksquare$
\begin{Remarque} The controls (\ref{EstimationLpprimeMolecules0}) and  (\ref{EstimationLpprimeMolecules1}) reflect a $\gamma$-h\"olderian gain of regularity of the solutions $\theta(t,\cdot)$ of the equation (\ref{EquationPrincipaleIntro}), however this gain is not instantaneous and some time $t> 2$ is needed to obtain the wished result. 
\end{Remarque}
We will prove now Theorem \ref{TheoLp'control} in two steps. First, in Section \ref{Secc_EvolMol} we study the evolution of the profile of the solutions $\psi$ to the dual equation (\ref{Evolution01}) and then in Section \ref{Secc_Iteration}  by a suitable iteration we will obtain the uniform bound (\ref{EstimationLpprimeMolecules1}).
\subsection{Molecule's evolution}\label{Secc_EvolMol}
The following theorem shows how the molecular properties are deformed with the evolution of the dual/linear equation (\ref{Evolution01}) for a time $0<s_0\ll 1$.
\begin{Theoreme}\label{SmallGeneralisacion} Assume that the hypotheses of Theorem \ref{THM_HOLDER} for the nonlinear equation (\ref{EquationPrincipaleIntro}) are in force. 
Consider $\psi_0$ a small molecule in the sense of Definition \ref{DefMolecules} for the local Hardy space $h^\mathfrak{s}(\mathbb{R}^n)$ where $\frac{n}{n+1}<\mathfrak{s}<1$ and let $\psi(s_0, x)$ be a solution at time $s_0$ of the dual problem (\ref{Evolution01}) associated with this molecular initial data $\psi_0$.\\

\noindent Then there exist positive constants $K$ and $\mathfrak{e}$ small enough such that for all time $0< s_0\le \mathfrak{e} r^\alpha$, we have the following estimates:
\begin{eqnarray}
\int_{\mathbb{R}^n}|\psi(s_0,x)||x-x(s_0)|^\omega dx &\leq &\big( (\zeta r)^\alpha+Ks_0\big)^{\frac{\omega-\gamma}{\alpha}}  \label{SmallConcentration},\\[4mm]
\|\psi(s_0,\cdot)\|_{L^\infty}&\leq & \frac{1}{\big((\zeta r)^\alpha+Ks_0\big)^{\frac{n+\gamma}{\alpha}}}\label{SmallLinftyevolution},\\
\|\psi(s_0,\cdot)\|_{L^1} &\leq & \frac{2v_n^{\frac{\omega}{n+\omega}}}{\big( (\zeta r)^\alpha+Ks_0\big)^{\frac{\gamma}{ \alpha}}},\label{SmallL1evolution}
\end{eqnarray}
where $\gamma$ is defined in (\ref{FormuleGammaRegulariteHolderHardy}), $\omega$ is the technical parameter given in Definition \ref{DefMolecules}, $\alpha=1+\varepsilon$ is the smoothness degree of the diffusion operator $\mathcal{L}^{\alpha}$ and $v_n$ denotes the volume of the $n$-dimensional unit ball.\\ 

\noindent The new molecule's center $x(s_0)$ used in formula (\ref{SmallConcentration}) is given by the evolution of the differential system
\begin{equation}\label{Defpointx_0Nonlinear}
\begin{cases}
x'(s)&= \bar{\mathbb{A}}_{[\theta]}(t-s,x(s))=\displaystyle{\int_{\mathbb{R}^{n}}}\mathbb{A}_{[\theta]}(t-s, x(s)-y)\varphi_{\rho_{0}}(y)dy,\quad  s\in [0,s_0],\\[2mm]
x(0)&= x_0,
\end{cases}
\end{equation}
where $\varphi_{\rho_{0}}(x)=\frac{1}{\rho_{0}^{n}}\varphi\big(\textcolor{black}{\frac x \rho_{0}}\big)$ (recall that $\mathbb{A}_{[\theta]}$ is given in (\ref{DefinitionNoyau})). Here $\varphi$ is a positive function in $\mathcal{C}^{\infty}_{0}(\mathbb{R}^{n})$ such that $supp(\varphi)\subset B(0,1)$ with $\rho_{0}=\zeta r\ll1$.
\end{Theoreme}
\textbf{Proof of the Theorem \ref{SmallGeneralisacion}.}
We will adopt the following strategy: we study first the Concentration condition (\ref{SmallConcentration}) in Section \ref{Secc_ConditionConcentration} and then we prove the Height condition (\ref{SmallLinftyevolution}) in Section \ref{Secc_ConditiondeHauteur}. With these two conditions at hand, the $L^1$ estimate (\ref{SmallL1evolution}) will be easily obtained in Section \ref{Secc_EstimationL1}. 
\subsubsection{Concentration condition}\label{Secc_ConditionConcentration}
To establish \eqref{SmallConcentration}, we introduce for $s\in [0,s_0]$ the function $\Omega_{s}(x)=|x-x(s)|^{\omega}$. For a given molecular initial data $\psi_{0}$ we write $\psi_{0}(x)=\psi_{0,+}(x)-\psi_{0,-}(x)$ where the functions $\psi_{0,\pm}(x)\geq 0$ have disjoint support and we denote by $\psi_\pm(s_0,x)$ two solutions of the dual equation (\ref{Evolution01}) at time $s_0$ with $\psi_{\pm}(0,x)=\psi_{0,\pm}(x)$.\\
 
The starting point for our study is the following: for all $s\in[0,s_0] $ we have
$$\int_{\mathbb{R}^n}|\psi(s_0,x)|\Omega_{s_0}(x)dx=\int_{\mathbb{R}^n}|\psi_{0}(x)|\Omega_{0}(x)dx+\int_0^{s_0} \partial_s\left(\int_{\R^n}|\psi(s,x)|\Omega_{s}(x) dx\right) ds,$$
now by the positivity principle stated in Proposition \ref{PropoLineairePositivite} we have $\psi_\pm(s_0,x)\geq 0$, then
$$|\psi(s_0,x)|=|\psi_+(s_0,x)-\psi_-(s_0,x)|\leq \psi_+(s_0,x)+\psi_-(s_0,x),$$ 
and we can thus write 
\begin{eqnarray*}
\int_{\mathbb{R}^n}|\psi(s_0,x)|\Omega_{s_0}(x)dx&\leq &\int_{\mathbb{R}^n}|\psi_{0}(x)|\Omega_{0}(x)dx\\
&&+\int_0^{s_0}\left|\partial_{s} \int_{\mathbb{R}^n}\psi_{+}(s,x)\Omega_{s}(x)dx\right|ds+\int_0^{s_0}\left|\partial_{s} \int_{\mathbb{R}^n}\psi_{-}(s,x)\Omega_{s}(x)dx\right|ds.
\end{eqnarray*}
We may assume, without loss of generality, that we have 
\begin{equation}\label{PREAL_ITE_GLOB_EN_TEMPS1}
\int_0^{s_0}\left|\partial_{s} \int_{\mathbb{R}^n}\psi_{-}(s,x)\Omega_{s}(x)dx\right|ds\leq \int_0^{s_0}\left|\partial_{s} \int_{\mathbb{R}^n}\psi_{+}(s,x)\Omega_{s}(x)dx\right|ds,
\end{equation}
and for the rest of the proof we will focus on the last term above. Since $\psi_{+}(s,\cdot)$ is the solution of the corresponding dual/linear equation (\ref{Evolution01}) with initial data $\psi_{0,+}(x)$,  for $s\in [0,s_0]$ we can write
\begin{eqnarray*}
I_s&:=&\left|\partial_{s} \int_{\mathbb{R}^n}\psi_{+}(s,x)\Omega_{s}(x)dx\right|\\
&=&\left|\int_{\mathbb{R}^n}\bigg(-\nabla\cdot\big(\mathbb{A}_{[\theta]}(t-s,x)\, \psi_{+}(s,x)\big)-\mathcal{L}^{\alpha}\psi_{+}(s,x)\bigg)\Omega_{s}(x)+\psi_{+}(s,x)\partial_{s} \Omega_{s}(x)dx\right|,
\end{eqnarray*}
recalling that $\Omega_{s}(x)=|x-x(s)|^{\omega}$ and reorganizing the previous terms we have
$$=\left|\int_{\mathbb{R}^n}-\nabla\Omega_{s}(x)\cdot x'(s)\psi_{+}(s,x)-\Omega_{s}(x)\nabla\cdot\big(\mathbb{A}_{[\theta]}(t-s,x)\, \psi_{+}(s,x)\big)-\Omega_{s}(x)\mathcal{L}^{\alpha}\psi_{+}(s,x)dx\right|.$$
By an integration by parts in the second term above and since the operator $\mathcal{L}^{\alpha}$ is symmetric we have
\begin{eqnarray*} 
I_s&=&\left|\int_{\mathbb{R}^n}-\nabla\Omega_{s}(x)\cdot x'(s)\psi_{+}(s,x)+\nabla\Omega_{s}(x) \cdot \mathbb{A}_{[\theta]}(t-s,x)\psi_{+}(s,x)-\mathcal{L}^{\alpha}\Omega_{s}(x)\psi_{+}(s,x)dx\right|\\
&\leq &\left|\int_{\mathbb{R}^n}\nabla\Omega_{s}(x)\cdot\bigg(\mathbb{A}_{[\theta]}(t-s,x)-x'(s)\bigg)\psi_{+}(s,x)dx\right|+\left|\int_{\mathbb{R}^n}\mathcal{L}^{\alpha}\Omega_{s}(x)\psi_{+}(s,x)dx\right|\\
&\leq& I_{s,1}+I_{s,2},
\end{eqnarray*}
and we obtain 
\begin{equation}\label{DECOUP_IS}
\int_{0}^{s_{0}}I_sds\leq \int_{0}^{s_{0}}I_{s,1}ds+\int_{0}^{s_{0}}I_{s,2}ds.
\end{equation}
We will study separately each of the quantities  $\displaystyle{\int_{0}^{s_{0}}I_{s,1}ds}$ and $\displaystyle{\int_{0}^{s_{0}}I_{s,2}ds}$ with two different propositions. The first contribution in \eqref{DECOUP_IS} is studied in Proposition \ref{Propo_Is1} below and constitutes the most technical part of the article since it requires to exploit in a very specific way the Besov stability controls stated in Corollary \ref{CorollaireLinftyBesov}. The contribution of the second term of  \eqref{DECOUP_IS} is easier to handle as it relies essentially on the properties of the operator $\mathcal{L}^{\alpha}$. It will be treated in Proposition \ref{Propo_Is2}.
\begin{Proposition}\label{Propo_Is1}
For the term $\displaystyle{\int_0^{s_0} I_{s,1}ds}$, under conditions (\ref{COND_B}) and (\ref{COND_E}) stated below,  we have
\begin{equation}\label{INTR_IS1_PROP}
\int_0^{s_0} I_{s,1}ds\le\Theta_1\; r^{\omega-\gamma-\alpha}\times s_0,
\end{equation}
where $\Theta_1$ is a small constant.
\end{Proposition}
{\bf Proof of Proposition \ref{Propo_Is1}.} 
With $\rho_0=\zeta r\ll 1$, we consider the following dyadic decomposition $\mathbb{R}^{n}=\displaystyle{B_{\rho_0}\cup \bigcup_{k\geq 1}E_{k}}$ where
\begin{equation}\label{DecoupageRn}
\begin{split}
B_{\rho_0}&=\{x\in \mathbb{R}^n: |x-x(s)|\leq \rho_0\} \quad \mbox{and}\\[2mm]
E_{k}&= \{x\in \mathbb{R}^n: 2^{k-1}\rho_0<|x-x(s)|\leq  2^{k}\rho_0\},\qquad \mbox{for } k\geq 1.
\end{split}
\end{equation}
Then we write
\begin{eqnarray}
I_{s,1}&=&\left|\int_{\mathbb{R}^n}\nabla\Omega_{s}(x)\cdot\bigg(\mathbb{A}_{[\theta]}(t-s,x)-x'(s)\bigg)\psi_{+}(s,x)dx\right|\nonumber\\
&=&\left|\int_{\mathbb{R}^n}\left(\mathds{1}_{B_{\rho_0}}+\sum_{k\geq 1}\mathds{1}_{E_{k}}\right)\nabla\Omega_{s}(x)\cdot\bigg(\mathbb{A}_{[\theta]}(t-s,x)-x'(s)\bigg)\psi_{+}(s,x)dx\right|\nonumber\\[2mm]
&\leq &\underbrace{\left|\int_{\mathbb{R}^n}\left(\mathds{1}_{B_{\rho_0}}\nabla\Omega_{s}(x)\cdot\big(\mathbb{A}_{[\theta]}(t-s,x)-\textcolor{black}{\bar{\mathbb{A}}_{[\theta]}(t-s,\cdot)}\big)\right)\psi_{+}(s,x)dx\right|}_{I_{s,B_{\rho_0}}}\label{EstimationBoule}\\
&&+\underbrace{\left|\int_{\mathbb{R}^n}\left(\sum_{k\geq 1}\mathds{1}_{E_{k,s}}\nabla\Omega_{s}(x)\cdot\big(\mathbb{A}_{[\theta]}(t-s,x)-\textcolor{black}{\bar{\mathbb{A}}_{[\theta]}(t-s,\cdot)}\big)\right)\psi_{+}(s,x)dx\right|}_{I_{s,E}},\label{EstimationCouronnes}
\end{eqnarray}
\textcolor{black}{recalling \eqref{Defpointx_0Nonlinear} for the last inequality}.
Thus we obtain
\begin{equation}\label{EstimationIntegreeTemps1} 
\int_0^{s_0} I_{s,1}ds\leq \int_0^{s_0} I_{s,B_{\rho_0}}ds+\int_0^{s_0} I_{s,E}ds.
\end{equation}
The analysis of the terms appearing in \eqref{EstimationBoule} and  \eqref{EstimationCouronnes} and represented in the previous right hand side is rather similar.
\textcolor{black}{The first key argument to handle those contributions is provided by Lemma \ref{LEMME_Lebesgue_Besov}, which precisely allows to control, on the elements of the previous spatial partition, a Lebesgue norm involving the difference of the drifts by a scaled Besov norm of the drift (which is again stable in this norm, see \eqref{HypoBesovBorne01})}. The controls obtained in Lemma \ref{LEMME_Lebesgue_Besov} are then used in Lemmas \ref{LEMMA_INTEGRANDS1} and \ref{LEMMA_INTEGRANDS3} to consider specifically the spatial domains appearing respectively in \eqref{EstimationBoule} and \eqref{EstimationCouronnes}.
\begin{Lemme}[From Lebesgue to Besov]\label{LEMME_Lebesgue_Besov}Let $1< p<+\infty$ and assume that the nonlinear drift $\mathbb{A}_{[\theta]}$ satisfies the hypothesis \A{B}-\A{C}. If moreover the term $\bar{\mathbb{A}}_{[\theta]}$ is defined as in equation (\ref{Defpointx_0Nonlinear}), then we have the following estimates over the sets $B_{\rho_{0}}$ and $E_{k}$, with $k\geq 1$, defined in expression (\ref{DecoupageRn}):
\begin{eqnarray}
\left\| \mathbb{A}_{[\theta]}(t-s,\cdot)-\bar{\mathbb{A}}_{[\theta]}(t-s,\cdot)\right\|_{L^{p}(B_{\rho_{0}})}&\leq &C\rho_{0}^{\frac{\sigma}{p}} \|\mathbb{A}_{[\theta]}(t-s,\cdot)\|_{\dot{B}^{\frac{\sigma}{p},p}_{p}}\label{EstimationLpBesovBoule}\\[3mm]
\|\mathbb{A}_{[\theta]}(t-s,\cdot)-\bar{\mathbb{A}}_{[\theta]}(t-s,\cdot)\|_{L^{p}(E_{k})}&\leq &C(2^{k}\rho_{0})^{\frac{\sigma}{p}}(2^{k})^{\frac{n}{p}}\|\mathbb{A}_{[\theta]}(t-s,\cdot)\|_{\dot{B}_{p}^{\frac{\sigma}{p},p}},\label{EstimationLpBesovGrandeCouronnes}
\end{eqnarray}
for some regularity index $0<\sigma$.
\end{Lemme}
The proof of this result is given in Lemma 4.1 of \cite{DCHSM_POTA}.\\
With this previous lemma, we are going to estimate each term of (\ref{EstimationIntegreeTemps1}). Although the general treatment of these terms is quite similar, some particularities must be taken into account. For the ball $B_{\rho_{1}}$ we have the following result.
\begin{Lemme}[Controls for the integrands over the ball $B_{\rho_{0}}$]\label{LEMMA_INTEGRANDS1}
With the notation of (\ref{DecoupageRn}) and \eqref{EstimationBoule} the following control holds for the small ball $B_{\rho_0}$:
\begin{eqnarray}\label{CTR_B}
\int_0^{s_0} I_{s,B_{\rho_0}} ds&\leq& \Theta_{1,1}\;r^{\omega-\gamma-\alpha}\times s_0,
\end{eqnarray}
where the technical quantity $\Theta_{1,1}=\Theta_{1,1}(\mu,\zeta,r)$ defined in the expression (\ref{DefTheta11}) below can be made small.
\end{Lemme}
{\bf Proof of Lemma \ref{LEMMA_INTEGRANDS1}.} We have
$$I_{s,B_{\rho_0}}=\left|\int_{\mathbb{R}^{n}}\left(\mathds{1}_{B_{\rho_0}}\nabla\Omega_{s}(x)\cdot\big(\mathbb{A}_{[\theta]}(t-s,x)-
\textcolor{black}{\bar{\mathbb{A}}_{[\theta]}(t-s,x)}
\big)\right)\psi_{+}(s,x)dx\right|,$$
and for $1<q_0,q_0'<+\infty$ with $\frac{1}{q_0}+\frac{1}{q_0'}=1$ we have
$$I_{s,B_{\rho_0}}\le \left\|\mathds{1}_{B_{\rho_0}} \nabla\Omega_{s}(\cdot)\cdot\big(\mathbb{A}_{[\theta]}(t-s,\cdot)-\textcolor{black}{\bar{\mathbb{A}}_{[\theta]}(t-s,\cdot)}
\big)\right\|_{L^{q_0'}}\|\psi_{+}(s,\cdot)\|_{L^{q_0}},$$
but since $|\nabla\Omega_{s}(x)|=c|x-x(s)|^{\omega-1}$
, we can write
$$I_{s,B_{\rho_0}}\leq C\left\||\cdot-x(s)|^{\omega-1}\big|\mathbb{A}_{[\theta]}(t-s,\cdot)-\bar{\mathbb{A}}_{[\theta]}(t-s,\cdot)\big|\right\|_{L^{q_0'}\textcolor{black}{(B_{\rho_0})}}\|\psi_{+}(s,\cdot)\|_{L^{q_0}}.$$
We need now to control the two previous Lebesgue norms and we need to handle the term $|x-x(s)|^{\omega-1}$ which is integrable over the ball $B_{\rho_{0}}$ under some conditions. Indeed, taking
\begin{equation}\label{ConditionBoule1}
\frac{1}{q_{0}'}=\frac{1}{a_{0}}+\frac{1}{p_0} \qquad \mbox{where}\qquad 1<a_{0}<\frac{n}{1-\omega}\quad  \mbox{and}\quad 2\leq p_0<+\infty,
\end{equation}
we thus derive from the H\"older inequality
$$I_{s,B_{\rho_0}}\leq C \||\cdot-x(s)|^{\omega-1}\|_{L^{a_0}(B_{\rho_0})}\big\|\mathbb{A}_{[\theta]}(t-s,\cdot)-\bar{\mathbb{A}}_{[\theta]}(t-s,\cdot)\big\|_{L^{p_0}(B_{\rho_0})}\|\psi_{+}(s,\cdot)\|_{L^{q_0}},$$
and after integration we obtain
$$I_{s,B_{\rho_0}}\leq C \rho_0^{\omega-1+\frac{n}{a_0}}\big\|\mathbb{A}_{[\theta]}(t-s,\cdot)-\bar{\mathbb{A}}_{[\theta]}(t-s,\cdot)\big\|_{L^{p_0}(B_{\rho_0})}\|\psi_{+}(s,\cdot)\|_{L^{q_0}}.$$
At this point, we apply the estimate (\ref{EstimationLpBesovBoule}) of Lemma \ref{LEMME_Lebesgue_Besov} with 
\begin{equation}\label{ConditionBoule2}
\frac{1}{p_0}=\frac{1-\nu_0}{\bf \bar p}+\frac{\nu_0}{2}\qquad\mbox{and}\qquad \sigma=p_0\nu_0(\frac{\alpha}{2}+\varepsilon_0)\quad \mbox{for some}\quad  0<\nu_0\leq 1,
\end{equation}
to obtain
\begin{eqnarray*}
I_{s,B_{\rho_0}}&\leq & C\rho_0^{\omega-1+\frac{n}{a_0}} \Big(\rho_{0}^{\nu_0(\frac{\alpha}{2}+\varepsilon_0)} \|\mathbb{A}_{[\theta]}(t-s,\cdot)\|_{\dot{B}_{p_0}^{\nu_0(\frac{\alpha}{2}+\varepsilon_0),p_0}}\Big)\|\psi_{+}(s,\cdot)\|_{L^{q_{0}}}\\
&\leq & C\rho_0^{\omega-1+\frac{n}{a_0}+\nu_0(\frac{\alpha}{2}+\varepsilon_0)} \|\mathbb{A}_{[\theta]}(t-s,\cdot)\|_{\dot{B}_{p_0}^{\nu_0(\frac{\alpha}{2}+\varepsilon_0),p_0}}\|\psi_{+}(s,\cdot)\|_{L^{q_{0}}}.
\end{eqnarray*}
Now, we integrate in the time variable this expression to obtain
$$\int_0^{s_0} I_{s,B_{\rho_0}}ds\leq C\rho_0^{\omega-1+\frac{n}{a_0}+\nu_0(\frac{\alpha}{2}+\varepsilon_0)} \int_0^{s_0} \|\mathbb{A}_{[\theta]}(t-s,\cdot)\|_{\dot{B}_{p_0}^{\nu_0(\frac{\alpha}{2}+\varepsilon_0),p_0}}\|\psi_{+}(s,\cdot)\|_{L^{q_{0}}}ds,$$
at this point we use the uniform in time estimates $\|\psi_{+}(s,\cdot)\|_{L^{q_{0}}}\leq \|\psi(s,\cdot)\|_{L^{q_{0}}}\leq \|\psi_0\|_{L^{q_{0}}}$ given by the maximum principle (\ref{ControlLqAprioriMolecules1}) and the estimates $\|\mathbb{A}_{[\theta]}(t-s,\cdot)\|_{\dot{B}_{p_0}^{\nu_0(\frac{\alpha}{2}+\varepsilon_0),p_0}}\leq C\mu^2$ given in (\ref{EstimationLinftyBesovPgrandA}) (recall that $2<t\leq N$ and $0\leq s\leq s_0$ with $0<s_0\ll1$) to write
\begin{eqnarray*}
\int_0^{s_0} I_{s,B_{\rho_0}}ds&\leq & C\rho_0^{\omega-1+\frac{n}{a_0}+\nu_0(\frac{\alpha}{2}+\varepsilon_0)}  \|\mathbb{A}_{[\theta]}\|_{L^\infty_s(\dot{B}_{p_0,x}^{\nu_0(\frac{\alpha}{2}+\varepsilon_0),p_0})}\|\psi_0\|_{L^{q_{0}}}\int_0^{s_0}ds\\
&\leq &C\rho_0^{\omega-1+\frac{n}{a_0}+\nu_0(\frac{\alpha}{2}+\varepsilon_0)}\mu^2\|\psi_{0}\|_{L^{q_{0}}} \times s_0.
\end{eqnarray*}
We recall now that $\rho_0=\zeta r\ll 1$, see (\ref{DecoupageRn}), and that $\|\psi_{0}\|_{L^{q_{0}}}\leq C(\zeta  r)^{-n+\frac{n}{q_0}-\gamma}$ by (\ref{Hipo4}) and we get 
$$\int_0^{s_0} I_{s,B_{\rho_0}}ds\leq C\mu^2(\zeta r)^{\omega-1+\frac{n}{a_0}+\nu_0(\frac{\alpha}{2}+\varepsilon_0)-n+\frac{n}{q_0}-\gamma}\times s_0.$$
Since $\frac{n}{a_0}-n+\frac{n}{q_0}=-\frac{n}{p_0}$ by (\ref{ConditionBoule1}), as we have $1-\frac{1}{q_0}=\frac{1}{q_0'}$, we can rewrite the previous inequality as follows
$$\int_0^{s_0} I_{s,B_{\rho_0}}ds\leq C\mu^2 (\zeta r)^{\frac{\nu_0\varepsilon_0}{2}}\times(\zeta r)^{\omega-\gamma-1-\frac{n}{p_0}+\nu_0\frac{\alpha+\varepsilon_0}{2}}\times s_0.$$
Now we define $\Theta_{1,1}$ by the expression 
\begin{equation}\label{DefTheta11}
\Theta_{1,1}=C\mu^2 (\zeta r)^{\frac{\nu_0\varepsilon_0}{2}}\zeta ^{\omega-\gamma-1-\frac{n}{p_0}+\nu_0\frac{\alpha+\varepsilon_0}{2}},
\end{equation}
and since $\rho_0=\zeta r\ll 1$ this quantity can be made very small and we have
$$\int_0^{s_0} I_{s,B_{\rho_0}}ds\leq \Theta_{1,1} r^{\omega-\gamma-1-\frac{n}{p_0}+\nu_0\frac{\alpha+\varepsilon_0}{2}}\times s_0,$$
thus, if we have the condition $\alpha\geq 1+\frac{n}{p_0}-\nu_0\frac{\alpha+\varepsilon_0}{2}$ which is equivalent to 
\begin{equation}\label{ConditionBoule3}
\alpha\geq \frac{p_0(2-\nu_0\varepsilon_0)+2n}{p_0(2+\nu_0)},
\end{equation}
we have $ r^{\omega-\gamma-1-\frac{n}{p_0}+\nu_0\frac{\alpha}{2}} \leq  r^{\omega-\gamma-\alpha}$ and we finally obtain
$$\int_0^{s_0} I_{s,B_{\rho_0}}ds\leq \Theta_{1,1} r^{\omega-\gamma-\alpha}\times s_0,$$
and the proof of Lemma \ref{LEMMA_INTEGRANDS1} is finished. \hfill $\blacksquare$
\begin{Remarque} In the proof of Lemma \ref{LEMMA_INTEGRANDS1} we used the following conditions on the parameters:
\begin{equation}\label{COND_B} 
\begin{cases}
\textcolor{black}{p_0>2,\ \bf\bar p>2, }\\[3mm]
1<q_{0}'\leq 2\leq  q_{0}<+\infty\quad \mbox{with }\quad \frac{1}{q_{0}}+\frac{1}{q_{0}'}=1,\\[3mm]
\frac{1}{q_{0}'}=\frac{1}{a_{0}}+\frac{1}{p_0} \quad \mbox{with}\quad 1<a_{0}<\frac{n}{1-\omega},\\[3mm]  
\frac{1}{p_0}=\frac{1-\nu_0}{\bf \bar p}+\frac{\nu_0}{2}\qquad\mbox{and}\qquad \sigma=p_0\nu_0(\frac{\alpha}{2}+\varepsilon_0)\quad \mbox{for some}\quad  0<\nu_0\leq 1,\\[3mm]
\alpha\geq \frac{p_0(2-\nu_0\varepsilon_0)+2n}{p_0(2+\nu_0)}.
\end{cases}
\tag{${\mathscr C}_B $}
\end{equation}
\end{Remarque}
\noindent We now study the quantity $\displaystyle{\int_0^{s_0}I_{s,E}ds}$ over the set $E=\displaystyle{\bigcup_{k\ge 1} E_k }$.
\begin{Lemme}[Controls for the integrands over the Dyadic coronas $\bigcup_{k\ge 1} E_k $]\label{LEMMA_INTEGRANDS3}
$$\int_0^{s_0}I_{s,E}ds\leq  \Theta_{1,2}\; r^{\omega-\gamma-\alpha}\times s_{0},$$
where the technical quantity $\Theta_{1,2}=\Theta_{1,2}(\mu,\zeta,r)$ defined in the expression (\ref{DefTheta12}) below can be made small.
\end{Lemme}
\textbf{Proof of Lemma \ref{LEMMA_INTEGRANDS3}.}  Since the set $E=\displaystyle{\bigcup_{k\ge 1} E_k }$ is the union of disjoint dyadic coronas (see (\ref{DecoupageRn})), we will need to derive estimates over each of the sets $(E_{k})_{k\geq 1}$ and thus to deal with some convergence issues. For this we write now, from the definition in \eqref{EstimationCouronnes}:
\begin{eqnarray*}
I_{s,E}&=&\left|\int_{\R^n}\left(\sum_{k\geq 1}\mathds{1}_{E_{k}}\nabla\Omega_{s}(x)\cdot\big(\mathbb{A}_{[\theta]}(t-s,x)-x'(s)\big)\right)\psi_{+}(s,x)dx\right|\\
&\le & \sum_{k\geq 1}\left|\int_{\R^n}\left(\mathds{1}_{E_{k}}\nabla\Omega_{s}(x)\cdot\big(\mathbb{A}_{[\theta]}(t-s,x)-x'(s)\big)\right)\psi_{+}(s,x)dx\right|.
\end{eqnarray*}
Now, by the H\"older inequality with $\frac{1}{q_1}+\frac{1}{q_1'}=1$ and $1<q_1<2<q_1'<+\infty$, we have
$$I_{s,E}\leq  \sum_{k\geq 1}\left\|\nabla\Omega_{s}(\cdot)\cdot\big(\mathbb{A}_{[\theta]}(t-s,\cdot)-x'(s)\big)\right\|_{L^{q_1'}(E_k)}\|\psi_{+}(s,\cdot)\|_{L^{q_1}},$$
and since $|\nabla\Omega_{s}(x)|=c|x-x(s)|^{\omega-1}$ and $x'(s)=\bar{\mathbb{A}}_{[\theta]}(t-s,x)$, we can write
$$I_{s,E}\leq  \sum_{k\geq 1}\left\||\cdot-x(s)|^{\omega-1}|\mathbb{A}_{[\theta]}(t-s,\cdot)-\bar{\mathbb{A}}_{[\theta]}(t-s,\cdot)|\right\|_{L^{q_1'}(E_k)}\|\psi_{+}(s,\cdot)\|_{L^{q_1}}.$$
Since by the definition (\ref{DecoupageRn}) of the dyadic corona $E_{k}$ we have $|x-x(s)|^{\omega-1}\leq C(2^{k}\rho_{0})^{\omega-1}$, we have
$$I_{s,E}\leq  C\sum_{k\geq 1}(2^{k}\rho_{0})^{\omega-1} \left\|\mathbb{A}_{[\theta]}(t-s,\cdot)-\bar{\mathbb{A}}_{[\theta]}(t-s,\cdot)\right\|_{L^{q_1'}(E_k)}\|\psi_{+}(s,\cdot)\|_{L^{q_1}},$$
and with the help of the control (\ref{EstimationLpBesovGrandeCouronnes}) given in Lemma \ref{LEMME_Lebesgue_Besov} with 
$$\frac{1}{q_1'}=\frac{1-\nu_1}{\bf \bar p}+\frac{\nu_1}{2}\qquad\mbox{and}\qquad \sigma=q_1'\nu_1(\frac{\alpha}{2}+\varepsilon_0)\quad \mbox{for some}\quad  0<\nu_1<1,$$
we can write
\begin{eqnarray*}
I_{s,E}&\leq & C\sum_{k\geq 1}(2^{k}\rho_{0})^{\omega-1} \Big( (2^{k}\rho_{0})^{\nu_1(\frac{\alpha}{2}+\varepsilon_0)}(2^{k})^{\frac{n}{q_{1}'}}\times \|\mathbb{A}_{[\theta]}(t-s,\cdot)\|_{\dot{B}^{\nu_1(\frac{\alpha}{2}+\varepsilon_0),q_{1}'}_{q_{1}'}}\Big)
\|\psi_{+}(s,\cdot)\|_{L^{q_1}}\\
&\leq & C\|\mathbb{A}_{[\theta]}(t-s,\cdot)\|_{\dot{B}^{\nu_1(\frac{\alpha}{2}+\varepsilon_0),q_{1}'}_{q_{1}'}}\|\psi_{+}(s,\cdot)\|_{L^{q_1}}\rho_{0}^{\omega-1+\nu_1(\frac{\alpha}{2}+\varepsilon_0)}\sum_{k\geq 1}2^{k\big(\omega-1+\nu_1(\frac{\alpha}{2}+\varepsilon_0)+\frac{n}{q_{1}'}\big)},
\end{eqnarray*}
and the condition
$$\omega-1+\nu_1(\frac{\alpha}{2}+\varepsilon_0)+\frac{n}{q_{1}'}<0,$$
guarantees that the previous sum converges. Now integrating with respect to the time variable we have
$$\int_0^{s_0}I_{s,E}ds\leq C\rho_{0}^{\omega-1+\nu_1(\frac{\alpha}{2}+\varepsilon_0)}\int_0^{s_0} \|\mathbb{A}_{[\theta]}(t-s,\cdot)\|_{\dot{B}^{\nu_1(\frac{\alpha}{2}+\varepsilon_0),q_{1}'}_{q_{1}'}}\|\psi_{+}(s,\cdot)\|_{L^{q_1}}ds.$$
Since we have $\|\psi_{+}(s,\cdot)\|_{L^{q_1}}\leq \|\psi(s,\cdot)\|_{L^{q_1}}\leq \|\psi_{0}\|_{L^{q_1}}$ by the maximum principle 
(\ref{ControlLqAprioriMolecules1}) and we have the estimates $\|\mathbb{A}_{[\theta]}(t-s,\cdot)\|_{\dot{B}_{q_1'}^{\nu_1(\frac{\alpha}{2}+\varepsilon_0),q_1'}}\leq C\mu^2$ given in (\ref{EstimationLinftyBesovPgrandA}) (recall that $t\in [1, N]$ and $0\leq s\leq s_0$ with $0<s_0\ll1$), we thus get
$$\int_0^{s_0}I_{s,E}ds\leq C\rho_{0}^{\omega-1+\nu_1(\frac{\alpha}{2}+\varepsilon_0)}\mu^2\|\psi_{0}\|_{L^{q_1}}\int_0^{s_0}ds=C\rho_{0}^{\omega-1+\nu_1(\frac{\alpha}{2}+\varepsilon_0)}\mu^2\|\psi_{0}\|_{L^{q_1}}\times s_0.$$
We recall now that $\rho_0=\zeta r\ll 1$ and that by (\ref{Hipo4}) we have $\|\psi_0\|_{L^{q_1}} \leq C (\zeta  r)^{-n+\frac{n}{q_1}-\gamma}$, so we can write
$$\int_0^{s_0}I_{s,E}ds\leq C\mu^2(\zeta  r)^{\omega-1+\nu_1(\frac{\alpha}{2}+\varepsilon_0)-n+\frac{n}{q_1}-\gamma}\times s_0.$$
Now we define $\Theta_{1,2}$ by the expression 
\begin{equation}\label{DefTheta12}
\Theta_{1,2}=C\mu^2 (\zeta r)^{\frac{\nu_1\varepsilon_0}{2}}\zeta ^{\omega-1+\nu_1\frac{\alpha+\varepsilon_0}{2}-n+\frac{n}{q_1}-\gamma},
\end{equation}
which can be made small since $\zeta r\ll 1$, and we  have 
$$\int_0^{s_0}I_{s,E}ds\leq \Theta_{1,2} r^{\omega-1+\nu_1(\frac{\alpha+\varepsilon_0}{2})-n+\frac{n}{q_1}-\gamma}\times s_0.$$
Note now that if  the condition
\begin{equation}\label{ConditionE3}
\alpha\geq \frac{q_1'(2-\nu_1\varepsilon_0)+2n}{q_1'(2+\nu_1)},
\end{equation}
holds, we then have $r^{\omega-1+\nu_1(\frac{\alpha+\varepsilon_0}{2})-n+\frac{n}{q_1}-\gamma}\leq  r^{\omega-\gamma-\alpha}$. We thus finally obtain
$$\int_0^{s_0}I_{s,E}ds\leq  \Theta_{1,2} r^{\omega-\gamma-\alpha}\times s_0,$$
and Lemma \ref{LEMMA_INTEGRANDS3} is proven. \hfill $\blacksquare$
\begin{Remarque} 

In the proof of Lemma \ref{LEMMA_INTEGRANDS3} we used the following conditions on the parameters:
\begin{equation}\label{COND_E} 
\begin{cases}
\textcolor{black}{\bf \bar p}>2,\\[3mm]
\frac{1}{q_1}+\frac{1}{q_1'}=1 \quad \mbox{and} \quad1<q_1<2<q_1'<+\infty,\\[3mm]
\frac{1}{q_1'}=\frac{1-\nu_1}{\bf \bar p}+\frac{\nu_1}{2}\qquad\mbox{and}\qquad \sigma=q_1'\nu_1(\frac{\alpha}{2}+\varepsilon_0)\quad \mbox{for some}\quad  0<\nu_1<1,\\[3mm]
\omega-1+\nu_1(\frac{\alpha}{2}+\varepsilon_0)+\frac{n}{q_{1}'}<0,\\[3mm]
\alpha\geq \frac{q_1'(2-\nu_1\varepsilon_0)+2n}{q_1'(2+\nu_1)}.
\end{cases}
\tag{${\mathscr C}_E $}
\end{equation}
\end{Remarque}
\textbf{End of the proof of Proposition \ref{Propo_Is1}.} From the estimates of Lemma \ref{LEMMA_INTEGRANDS1} and \ref{LEMMA_INTEGRANDS3}, and under the conditions stated in (\ref{COND_B}) and (\ref{COND_E}) for the involved parameters, we then derive:
$$\int_{0}^{s_{0}}I_{s,1}ds\leq \int_{0}^{s_{0}}I_{s,B}+I_{s,E}\,ds\leq \left(\Theta_{1,1}+\Theta_{1,2}\right)\times r^{\omega-\gamma-\alpha} \times s_0,$$
and setting $\Theta_1:=\Theta_{1,1}+\Theta_{1,2}\ll1$, this completes the proof of Proposition \ref{Propo_Is1}. \hfill $\blacksquare$\\

The contribution $I_{s,2}$ in (\ref{DECOUP_IS}) only depends on the operator. Precisely, we get:
\begin{Proposition}\label{Propo_Is2}
For the second integral $\displaystyle{\int_{0}^{s_{0}}I_{s,2}ds}$ in (\ref{DECOUP_IS}), we have:
$$\int_{0}^{s_{0}}I_{s,2}ds\leq  \Theta_2\times r^{\omega-\gamma-\alpha}\times s_{0},$$
where $\Theta_2=\Theta_2(\zeta)$ is given by
$$\Theta_2(\zeta)=C\left( \zeta^{\mathfrak{m}_{0}(\omega-\alpha+{n})-n-\gamma}+\zeta^{\mathfrak{m}_{0}(1+\epsilon)[\omega-\alpha-\frac{n}{\mathfrak{p}_{1}}]-n+\frac{n}{\mathfrak{p}'_{1}}-\gamma} +\zeta^{\mathfrak{m}_{1}[\omega-\alpha+\frac{n}{\mathfrak{p}_{2}}]-n+\frac{n}{\mathfrak{p}'_{2}}-\gamma}\right),$$
and the parameters above satisfy 
\begin{equation}\label{Cond_Generales2}
\begin{cases}
0<\mathfrak{m}_{0}<1<\mathfrak{m}_{1}, \,1<\mathfrak{p}_{1},\mathfrak{p}_{1}',\mathfrak{p}_{2},\mathfrak{p}_{2}'<+\infty \quad\mbox{and}\quad \frac{1}{\mathfrak{p}_{1}}+ \frac{1}{\mathfrak{p}_{1}'}= \frac{1}{\mathfrak{p}_{2}}+ \frac{1}{\mathfrak{p}_{2}'}=1,\\[3mm]
\frac{n}{\alpha-\omega}<\mathfrak{p}_{1}, \mathfrak{p}_{2}<+\infty \quad\mbox{and}\quad\epsilon=\frac{\ln(1-\zeta^{(\mathfrak{m}_{1}-\mathfrak{m}_{0})(\mathfrak{p}_{1}(\omega-\alpha)+n)})}{(\mathfrak{p}_{1}(\omega-\alpha)+n)\mathfrak{m}_{0}\ln(\zeta)}, 
\end{cases}
\tag{${\mathscr C}_2$}
\end{equation}
and $\Theta_2$ is such that the quantity $\frac{\Theta_2(\zeta)}{\zeta^{\omega-\gamma-\alpha}}$ can be made very small for $\zeta $ large enough.
\end{Proposition}
{\bf Proof of Proposition \ref{Propo_Is2}.}  It is enough to follow essentially the same ideas used for Proposition \ref{Propo_Is1} with minor modifications which make computations much easier: indeed, for $\zeta\gg 1$, $0<r<1$ and $0<\mathfrak{m}_{0}<1<\mathfrak{m}_{1}$, we define $\mathfrak{r}_{0}, \mathfrak{r}_{1}$ by the expressions
\begin{equation}\label{DefinitionRhoMIs2}
\begin{split}
\mathfrak{r}_{0}=\zeta^{\mathfrak{m}_{0}}r\leq 1,
\mathfrak{r}_{1}=\zeta^{\mathfrak{m}_{1}}r\leq 1,
\end{split}
\end{equation}
and we consider the decomposition $\displaystyle{\mathbb{R}^{n}=B_{\mathfrak{r}_{0}}\cup\mathcal{C}_{(\mathfrak{r}_{0}, \mathfrak{r}_{1})}\cup \big(\bigcup_{k\geq 1}E_{k}\big)}$ where $B_{\mathfrak{r}_{0}}, \mathcal{C}_{(\mathfrak{r}_{0}, \mathfrak{r}_{1})}$ and $E_{k}$ are defined  in a very similar way as in (\ref{DecoupageRn}) but replacing $\rho_{0}, \rho_{1}$ by $\mathfrak{r}_{0}, \mathfrak{r}_{1}$. Indeed we have
\begin{equation*}
\begin{split}
B_{\mathfrak{r}_{0}}&=\{x\in \mathbb{R}^n: |x-x(s)|\leq \mathfrak{r}_0\}\\[2mm]
\mathcal{C}_{(\mathfrak{r}_{0}, \mathfrak{r}_{1})}&=\{x\in \mathbb{R}^n: \mathfrak{r}_0<|x-x(s)|\leq \mathfrak{r}_1\} \quad \mbox{and}\\[2mm]
E_{k}&= \{x\in \mathbb{R}^n: 2^{k-1} \mathfrak{r}_1<|x-x(s)|\leq  2^{k} \mathfrak{r}_1\},\qquad \mbox{for } k\geq 1.
\end{split}
\end{equation*}
We thus can write
\begin{eqnarray}
I_{s,2}&=&\left|\int_{\mathbb{R}^n}\mathcal{L}^{\alpha}\Omega_{s}(x)\psi_{+}(s,x)dx\right|\notag\\
&\leq &\left|\int_{\mathbb{R}^n}\mathds{1}_{B_{\mathfrak{r}_{0}}}\mathcal{L}^{\alpha}\Omega_{s}(x)\psi_{+}(s,x)dx\right|+\left|\int_{\mathbb{R}^n}\mathds{1}_{\mathcal{C}_{(\mathfrak{r}_{0}, \mathfrak{r}_{1})}}\mathcal{L}^{\alpha}\Omega_{s}(x)\psi_{+}(s,x)dx\right|\notag\\
& &+\left|\int_{\mathbb{R}^n}\sum_{k\geq 1}\mathds{1}_{E_{k}}\mathcal{L}^{\alpha}\Omega_{s}(x)\psi_{+}(s,x)dx\right|.\label{EstimationBoulePetiteCouronneGrandeCouronnes}
\end{eqnarray}
\begin{itemize}
\item[$\bullet$] For the first term of (\ref{EstimationBoulePetiteCouronneGrandeCouronnes}) we write
$$\left|\int_{\mathbb{R}^n}\mathds{1}_{B_{\mathfrak{r}_{0}}}\mathcal{L}^{\alpha}\Omega_{s}(x)\psi_{+}(s,x)dx\right|\leq \|\mathcal{L}^{\alpha}|\cdot-x(s)|^{\omega}\|_{L^{1}(B_{\mathfrak{r}_0})}\|\psi_{+}(s,\cdot)\|_{L^{\infty}},
$$
note that the last term above is equivalent up to a change of variables to
$\displaystyle{\int_{B(0,\mathfrak{r}_0)} |\mathcal{L}^{\alpha}|x|^{\omega}|dx}$, and by classical homogeneity arguments (see the book \cite{Grafakos}) we have $\big|\mathcal{L}^{\alpha}|x|^{\omega}\big|=C|x|^{\omega-\alpha}$ and we obtain $
\|\mathcal{L}^{\alpha}|x-x(s)|^{\omega}\|_{L^{1}(B_{\mathfrak{r}_0})}\leq C \mathfrak{r}_0^{\omega-\alpha+{n}}$. Thus by the maximum principle, by the molecular hypothesis (\ref{Hipo2}) and by the definition of $\mathfrak{r}_{0}$ given in  (\ref{DefinitionRhoMIs2}) we obtain
\begin{equation}\label{EstimationISB2AvantIntegrationTemps}
\left|\int_{\mathbb{R}^n}\mathds{1}_{B_{\mathfrak{r}_{0}}}\mathcal{L}^{\alpha}\Omega_{s}(x)\psi_{+}(s,x)dx\right|\leq C \zeta^{\mathfrak{m}_{0}(\omega-\alpha+{n})-n-\gamma}\times r^{\omega-\gamma-\alpha}.
\end{equation}
\item[$\bullet$]  For the second term of (\ref{EstimationBoulePetiteCouronneGrandeCouronnes}), we set $1<\mathfrak{p}_{1}, \mathfrak{p}'_{1}<+\infty$ such that $\frac{1}{\mathfrak{p}_{1}}+\frac{1}{\mathfrak{p}'_{1}}=1$ and $\frac{n}{\alpha-\omega}<\mathfrak{p}_{1}<+\infty$, and we write
\begin{eqnarray*}
\left|\int_{\mathbb{R}^n}\mathds{1}_{\mathcal{C}_{(\mathfrak{r}_{0}, \mathfrak{r}_{1})}}\mathcal{L}^{\alpha}\Omega_{s}(x)\psi_{+}(s,x)dx\right|&\leq &\|\mathcal{L}^{\alpha}\Omega_{s}(x)\|_{L^{\mathfrak{p}_{1}}(\mathcal{C}_{(\mathfrak{r}_{0}, \mathfrak{r}_{1})})} \|\psi_{+}(s\cdot)\|_{L^{\mathfrak{p}'_{1}}}\\
&\leq &\|\mathcal{L}^{\alpha}\Omega_{s}(x)\|_{L^{\mathfrak{p}_{1}}(\mathcal{C}_{(\mathfrak{r}_{0}, \mathfrak{r}_{1})})} \|\psi_{0}\|_{L^{\mathfrak{p}'_{1}}},
\end{eqnarray*}
where we used the maximum principle for the function $\psi_{+}(s,\cdot)$. Now, by homogeneity we have $\|\mathcal{L}^{\alpha}\Omega_{s}(x)\|_{L^{\mathfrak{p}_{1}}(\mathcal{C}_{(\mathfrak{r}_{0}, \mathfrak{r}_{1})})}=C\||\cdot-x(s)|^{\omega-\alpha}\|_{L^{\mathfrak{p}_{1}}(\mathcal{C}_{(\mathfrak{r}_{0}, \mathfrak{r}_{1})})}\leq C(\mathfrak{r}_{0}^{\mathfrak{p}_1(\omega-\alpha)+n}-\mathfrak{r}_{1}^{\mathfrak{p}_{1}(\omega-\alpha)+n})^{\frac{1}{\mathfrak{p}_{1}}}$ and if we consider
$$\epsilon=\frac{\ln(1-\zeta^{(\mathfrak{m}_{1}-\mathfrak{m}_{0})(\mathfrak{p}_{1}(\omega-\alpha)+n)})}{(\mathfrak{p}_{1}(\omega-\alpha)+n)\mathfrak{m}_{0}\ln(\zeta)},$$ 
we have then by the definition of $\mathfrak{r}_{0}$ and $\mathfrak{r}_{1}$ given in (\ref{DefinitionRhoMIs2}):
$$(\mathfrak{r}_{0}^{\mathfrak{p}_{1}(\omega-\alpha)+n}-\mathfrak{r}_{1}^{\mathfrak{p}_{1}(\omega-\alpha)+n})^{\frac{1}{\mathfrak{p}_{1}}}=\zeta^{\mathfrak{m}_{0}(1+\epsilon)[\omega-\alpha-\frac{n}{\mathfrak{p}_{1}}]} r^{\omega-\alpha+\frac{n}{\mathfrak{p}_{1}}},$$
and thus using (\ref{Hipo4}) to estimate the molecular initial data $ \|\psi_{0}\|_{L^{\mathfrak{p}'_{1}}}$ we have
\begin{eqnarray}
\left|\int_{\mathbb{R}^n}\mathds{1}_{\mathcal{C}_{(\mathfrak{r}_{0}, \mathfrak{r}_{1})}}\mathcal{L}^{\alpha}\Omega_{s}(x)\psi_{+}(s,x)dx\right|&\leq & C \zeta^{\mathfrak{m}_{0}(1+\epsilon)[\omega-\alpha-\frac{n}{\mathfrak{p}_{1}}]} r^{\omega-\alpha+\frac{n}{\mathfrak{p}_{1}}}\times (\zeta r)^{-n+\frac{n}{\mathfrak{p}'_{1}}-\gamma}\notag\\
&\leq  &C\zeta^{\mathfrak{m}_{0}(1+\epsilon)[\omega-\alpha-\frac{n}{\mathfrak{p}_{1}}]-n+\frac{n}{\mathfrak{p}'_{1}}-\gamma} \times r^{\omega-\gamma-\alpha}.\label{EstimationISC2AvantIntegrationTemps}
\end{eqnarray}
\item[$\bullet$]  For the last term of (\ref{EstimationBoulePetiteCouronneGrandeCouronnes}), by the H\"older inequality with $1<\mathfrak{p}_{2}, \mathfrak{p}'_{2}<+\infty$ such that $\frac{1}{\mathfrak{p}_{2}}+\frac{1}{\mathfrak{p}'_{2}}=1$ and $\omega-\alpha+\frac{n}{\mathfrak{p}_{2}}<0$ and by homogeneity arguments, we have
\begin{eqnarray*}
\left|\int_{\mathbb{R}^n}\sum_{k\geq 1}\mathds{1}_{E_{k}}\mathcal{L}^{\alpha}\Omega_{s}(x)\psi_{+}(s,x)dx\right|&\leq &\sum_{k\geq 1}\||\cdot-x(s)|^{\omega-\alpha}\|_{L^{\mathfrak{p}_{2}}(E_{k})} \|\psi_{+}(s,\cdot)\|_{L^{\mathfrak{p}'_{2}}}\\
&\leq &C\sum_{k\geq 1} (2^{k}\mathfrak{r}_{1})^{\omega-\alpha+\frac{n}{\mathfrak{p}_{2}}} \|\psi_{+}(s,\cdot)\|_{L^{\mathfrak{p}'_{2}}},
\end{eqnarray*}
note that by the condition $\omega-\alpha+\frac{n}{\mathfrak{p}_{2}}$ the previous sum is convergent. Now, by the definition of $\mathfrak{r}_{1}$ given in (\ref{DefinitionRhoMIs2}), by the maximum principle for $\psi_{+}(s,\cdot)$ and by (\ref{Hipo4}) we obtain 
\begin{eqnarray}
\left|\int_{\mathbb{R}^n}\sum_{k\geq 1}\mathds{1}_{E_{k}}\mathcal{L}^{\alpha}\Omega_{s}(x)\psi_{+}(s,x)dx\right|&\leq &C\mathfrak{r}_{1}^{\omega-\alpha+\frac{n}{\mathfrak{p}_{2}}} \times (\zeta r)^{-n+\frac{n}{\mathfrak{p}'_{2}}-\gamma}\notag\\
&\leq &C\zeta^{\mathfrak{m}_{1}[\omega-\alpha+\frac{n}{\mathfrak{p}_{2}}]-n+\frac{n}{\mathfrak{p}'_{2}}-\gamma}\times r^{\omega-\gamma-\alpha}.  \label{EstimationISE2AvantIntegrationTemps}
\end{eqnarray}
\end{itemize}
Thus, with estimates (\ref{EstimationISB2AvantIntegrationTemps}), (\ref{EstimationISC2AvantIntegrationTemps}) and (\ref{EstimationISE2AvantIntegrationTemps}) we obtain the following control for (\ref{EstimationBoulePetiteCouronneGrandeCouronnes}):
$$I_{s,2}\leq C\left( \zeta^{\mathfrak{m}_{0}(\omega-\alpha+{n})-n-\gamma}+\zeta^{\mathfrak{m}_{0}(1+\epsilon)[\omega-\alpha-\frac{n}{\mathfrak{p}_{1}}]-n+\frac{n}{\mathfrak{p}'_{1}}-\gamma} +\zeta^{\mathfrak{m}_{1}[\omega-\alpha+\frac{n}{\mathfrak{p}_{2}}]-n+\frac{n}{\mathfrak{p}'_{2}}-\gamma}\right)\times r^{\omega-\gamma-\alpha},$$
which becomes after an integration in the time variable
$$\int_{0}^{s_{0}}I_{s,2}ds\leq C\left( \zeta^{\mathfrak{m}_{0}(\omega-\alpha+{n})-n-\gamma}+\zeta^{\mathfrak{m}_{0}(1+\epsilon)[\omega-\alpha-\frac{n}{\mathfrak{p}_{1}}]-n+\frac{n}{\mathfrak{p}'_{1}}-\gamma} +\zeta^{\mathfrak{m}_{1}[\omega-\alpha+\frac{n}{\mathfrak{p}_{2}}]-n+\frac{n}{\mathfrak{p}'_{2}}-\gamma}\right)\times r^{\omega-\gamma-\alpha}\times s_{0}.$$
Moreover, if we set
\begin{equation}\label{DefinitionTheta2Global}
\Theta_2(\zeta)=C\left( \zeta^{\mathfrak{m}_{0}(\omega-\alpha+{n})-n-\gamma}+\zeta^{\mathfrak{m}_{0}(1+\epsilon)[\omega-\alpha-\frac{n}{\mathfrak{p}_{1}}]-n+\frac{n}{\mathfrak{p}'_{1}}-\gamma} +\zeta^{\mathfrak{m}_{1}[\omega-\alpha+\frac{n}{\mathfrak{p}_{2}}]-n+\frac{n}{\mathfrak{p}'_{2}}-\gamma}\right),
\end{equation}
we finally obtain 
$$\int_{0}^{s_{0}}I_{s,2}ds\leq  \Theta_2\times r^{\omega-\gamma-\alpha}\times s_{0}.$$
To finish Proposition \ref{Propo_Is2} it only remains to show that the quantity 
\begin{equation}\label{FormuleAbsorptionFinale2}
\frac{\Theta_2(\zeta)}{\zeta^{\omega-\gamma-\alpha}},
\end{equation} 
can be made small if $\zeta$ is big enough: for this we only have to check the sign of the powers of $\zeta$ of this quantity, and from (\ref{DefinitionTheta2Global}) above we obtain the conditions
\begin{equation}\label{AbsConstantesIS2}
\begin{cases}
\mathfrak{m}_{0}(\omega-\alpha+{n})-n+\alpha-\omega<0,\\[3mm]
\mathfrak{m}_{0}(1+\epsilon)[\omega-\alpha-\frac{n}{\mathfrak{p}_{1}}]-n+\frac{n}{\mathfrak{p}'_{1}}+\alpha-\omega<0,\\[3mm]
\mathfrak{m}_{1}[\omega-\alpha+\frac{n}{\mathfrak{p}_{2}}]-n+\frac{n}{\mathfrak{p}'_{2}}+\alpha-\omega<0.
\end{cases}
\tag{${\mathscr{A}}$}
\end{equation}
\hfill $\blacksquare$
\begin{Remarque}\label{RemarqueRoleConditionIS2}
Not that, unlike Proposition \ref{Propo_Is1}, the constraints (\ref{AbsConstantesIS2}) above are essential to absorb the constants that appear in the computations.  
\end{Remarque}
We finish now the proof of the evolution of the concentration condition (\ref{SmallConcentration}) and from equation (\ref{PREAL_ITE_GLOB_EN_TEMPS1}) we write:
\begin{eqnarray*}
\int_{\mathbb{R}^n}|\psi(s_0,x)| \Omega_{s_{0}}(x) dx&\leq &\int_{\mathbb{R}^n}|\psi_{0}(x)|\Omega_{0}(x)dx+\int_0^{s_0}\left|\partial_{s} \int_{\mathbb{R}^n}\psi_{+}(s,x)\Omega_{s}(x)dx\right|ds\\
&&+\int_0^{s_0}\left|\partial_{s} \int_{\mathbb{R}^n}\psi_{-}(s,x)\Omega_{s}(x)dx\right|ds.
\end{eqnarray*}
Recall that $\displaystyle{\int_{\mathbb{R}^n}|\psi_{0}(x)|\Omega_{0}(x)dx=(\zeta r)^{\omega-\gamma}}$ (since the initial data is a molecule). \textcolor{black}{Since the time sensitivities of the positive and negative parts of $\psi $ are controlled in a similar way, we derive from Propositions \ref{Propo_Is1}, \ref{Propo_Is2} that setting} $\Theta:=\Theta_1+\Theta_2\ll 1$:  
\begin{eqnarray*}
\int_{\mathbb{R}^n}|\psi(s_0,x)| \Omega_{s_{0}}(x) dx&\le & (\zeta r)^{\omega-\gamma}+2\left(\Theta\times r^{\omega-\gamma-\alpha}\times s_0\right)\\
&\leq & (\zeta r)^{\omega-\gamma}\big(1+2\frac{\Theta}{\zeta^{\omega-\gamma}}\frac{s_0}{r^\alpha}\big).
\end{eqnarray*}
Recalling that $\Omega_{s_{0}}(x)=|x-x(s_0)|^{\omega}$ and that we assumed $0\leq s_0\leq \mathfrak{e} r^\alpha$ we can choose $ \mathfrak{e}$ small enough to have:
{\small
\begin{eqnarray}
\int_{{\mathbb R}^n}|\psi(s_0,x)| |x-x(s_0)|^\omega dx&\leq & (\zeta r)^{\omega-\gamma}\bigg(1+\textcolor{black}{4} \frac{\alpha}{\omega-\gamma}\frac{\Theta}{\zeta^{\omega-\gamma}} \frac{s_0}{r^\alpha}\bigg)^{\frac{\omega-\gamma}{\alpha}}\notag\\
&\leq & \bigg( (\zeta r)^\alpha+\textcolor{black}{4} \frac{\alpha}{\omega-\gamma}\frac{\Theta}{\zeta^{\omega-\gamma-\alpha}} s_0\bigg)^{\frac{\omega-\gamma}{\alpha}}\leq \big( (\zeta r)^\alpha+ K s_0\big)^{\frac{\omega-\gamma}{\alpha}},\label{PresqueK}
\end{eqnarray}
}
since we can set $\frac{\Theta}{\zeta^{\omega-\gamma-\alpha}}$ small enough in order to satisfy
\begin{equation}\label{ConditionK11}
\textcolor{black}{4} \frac{\alpha}{\omega-\gamma}\frac{\Theta}{\zeta^{\omega-\gamma-\alpha}}\leq K,
\end{equation} 
where $K$ in \eqref{PresqueK} is a (small) constant that intervenes in the height condition. Indeed, we have 
$$\frac{\Theta}{\zeta^{\omega-\gamma-\alpha}}=\frac{\Theta_1}{\zeta^{\omega-\gamma-\alpha}}+\frac{\Theta_2}{\zeta^{\omega-\gamma-\alpha}},$$
and we note that the term $\frac{\Theta_1}{\zeta^{\omega-\gamma-\alpha}}$ can be made very small as by (\ref{DefTheta11}), this quantity is governed by the parameter $0<r\ll 1$, moreover, by (\ref{FormuleAbsorptionFinale2}) we have that $\frac{\Theta_2}{\zeta^{\omega-\gamma-\alpha}}$ can also made very small and this concludes the proof of the concentration condition (\ref{SmallConcentration}) of Theorem \ref{SmallGeneralisacion}.\hfill $\blacksquare$
\subsubsection{Height condition}\label{Secc_ConditiondeHauteur}
We study in this section the Height condition (\ref{SmallLinftyevolution}) given by the expression
$$\|\psi(s_0,\cdot)\|_{L^\infty}\leq \frac{1}{\big((\zeta r)^\alpha+Ks_0\big)^{\frac{n+\gamma}{\alpha}}}.$$
To establish the control of the theorem we aim at proving that
\begin{equation}\label{INEG_DIFF}
\frac{d}{ds}\|\psi(s, \cdot)\|_{L^\infty}\leq -K\big(\frac{n+\gamma}{\alpha} \big)((\zeta r)^\alpha+Ks)^{-{(\omega-\gamma)}/{(n+\omega)}} \,\|\psi(s, \cdot)\|_{L^\infty}^{1+{\alpha}/{(n+\omega)}}.\!\!\!\!
\end{equation}
Indeed, integrating \eqref{INEG_DIFF} yields
\begin{eqnarray*}
\int_0^{s_0} \frac{d}{ds}\Big( \|\psi(s,\cdot)\|_{L^\infty}^{-{\alpha}/{(n+\omega)}}\Big)\,ds  &\geq &\int_0^{s_0}  \frac{d}{ds}\big( [(\zeta r)^\alpha+Ks]^{{(n+\gamma)}/{(n+\omega)}} \big)\,ds, 
\\
\|\psi(s_0,\cdot)\|_{L^\infty}^{-{\alpha}/{(n+\omega)}}&\geq&  [(\zeta r)^\alpha+Ks_0]^{(n+\gamma)/(n+\omega)} +\big(\|\psi(0,\cdot)\|_{L^\infty}^{-{\alpha}/{(n+\omega)}}-[(\zeta r)^\alpha]^{(n+\gamma)/(n+\omega)} \big)\\
&\geq & [(\zeta r)^\alpha+Ks_0]^{(n+\gamma)/(n+\omega)}.
\end{eqnarray*}
Recalling the initial height condition $\|\psi(0,\cdot)\|_{L^\infty}\leq (\zeta r)^{-(n+\gamma)} $ for the last inequality, we therefore derive
$$\|\psi(s_0,\cdot)\|_{L^\infty}\leq ((\zeta r)^\alpha+Ks_0)^{-(n+\gamma)/{\alpha}},$$
which is the required control.\\

Assume that the molecules we are working with are smooth enough and in particular continuous. Following an idea of \cite{Cordoba} (Section~4, p.~\mbox{522--523}) (see also \cite{Jacob}, p.\ 346), we will denote for $s\in [0,s_0] $ by $\overline{x}_s$ the point of $\mathbb{R}^n$ such that $\psi(s, \overline{x}_s)=\|\psi(s, \cdot)\|_{L^\infty}$. Thus we can write, by the properties (\ref{DefKernelLevy}) of the function $\pi$,
\begin{align}\nonumber
\frac{d}{ds}\|\psi(s, \cdot)\|_{L^\infty}&\leq -\int_{\mathbb{R}^n}[\psi(s, \overline{x}_s)-\psi(s, \overline{x}_s-y)]\, \pi(y)\,dy\\
&\leq -\mathfrak{c}_1\int_{\{|\overline{x}_s-y|<1\}}\frac{\psi(s, \overline{x}_s)-\psi(s,y)}{|\overline{x}_s-y|^{n+\alpha}}\,dy\leq 0.
\label{Infty1}
\end{align}
To establish the differential inequality \eqref{INEG_DIFF} for $s\in [0,s_0] $, let us first consider a corona centered in $\bar x_s $ defined by $\displaystyle{\mathcal{C}(R, \rho R)=\{y\in \mathbb{R}^n:R\leq|\overline{x}_s-y|\leq \rho R\}}$, where the parameter $R>0$ to be specified later on is such that $0<\rho R<1$ with $\rho>2$. Then,
\begin{equation*}
\int_{\{|\overline{x}_{s}-y|<1\}}\frac{\psi(s, \overline{x}_{s})-\psi(s, y)}{|\overline{x}_{s}-y|^{n+\alpha}}\,dy\geq \int_{\mathcal{C}(R,\rho R)}\frac{\psi(s,\overline{x}_{s})-\psi(s,y)}{|\overline{x}_{s}-y|^{n+\alpha}}\,dy.
\end{equation*}
Define now the sets $B_1$ and $B_2$ by
\begin{align*}
B_1&=\{y\in \mathcal{C}(R,\rho R): \psi(s,\overline{x}_{s})-\psi(s,y)\geq \tfrac{1}{2}\psi(s,\overline{x}_{s})\},
 \\[3mm]
B_2&=\{y\in \mathcal{C}(R,\rho R): \psi(s,\overline{x}_{s})-\psi(s,y)< \tfrac{1}{2}\psi(s,\overline{x}_{s})\},
\end{align*}
 such that $\mathcal{C}(R,\rho R)=B_1\cup B_2$. We obtain then the inequalities
\begin{align}\nonumber
\int_{\mathcal{C}(R,\rho R)}&\,\frac{\psi(s,\overline{x}_{s})-\psi(s, y)}{ |\overline{x}_{s}-y|^{n+\alpha}}\,dy \geq \int_{B_1}\frac{\psi(s,\overline{x}_{s})-\psi(s, y)}{|\overline{x}_{s}-y|^{n+\alpha}}\,dy \geq \frac{\psi(s,\overline{x}_{s})}{2\rho^{n+\alpha} R^{n+\alpha}}|B_1|
\\
&=\frac{\psi(s,\overline{x}_{s})}{2\rho^{n+\alpha} R^{n+\alpha}}\left(|\mathcal{C}(R,\rho R)|-|B_2|\right)
\geq \frac{\psi(s,\overline{x}_{s})}{2\rho^{n+\alpha}R^{n+\alpha}}\big(v_n(\rho^n {-}1)R^n-|B_2|\big),\label{Infty3}
\end{align}
where $v_n$ denotes the volume of the unit ball.\\[5mm]

To continue, we need to estimate the quantity $|B_2|$ in the right-hand side of~\eqref{Infty3} in terms of $\psi(s,\overline{x}_{s})$ and $R$. We will distinguish two cases and prove the following estimates:
\begin{enumerate}
\item[1)] If $|\overline{x}_{s}-x(s)|>2\rho R$ or $|\overline{x}_{s}-x(s)|<R/2$, then
\begin{equation}\label{FormEstima1}
C_1\big((\zeta r)^{\alpha}+K s\big)^{(\omega-\gamma)/{\alpha}}\, \psi(s,\overline{x}_{s})^{-1}\, R^{-\omega}\geq |B_2|.
\end{equation}
\item[2)] If $R/2\leq |\overline{x}_{s}-x(s)|\leq 2\rho R$, then
\begin{equation}\label{FormEstima2}
\big(C_2   \,\big((\zeta r)^{\alpha}+K s\big)^{(\omega-\gamma)/{\alpha}}\, R^{n-\omega}\, \psi(s,\overline{x}_{s})^{-1}\big)^{1/2}\geq |B_2|.
\end{equation}
\end{enumerate}

For these two controls, our starting point is the concentration condition established in Theorem \ref{SmallGeneralisacion}, eq. \eqref{SmallConcentration}; indeed we can write
\begin{align}\nonumber
\big((\zeta r)^{\alpha}+K s\big)^{(\omega-\gamma)/{\alpha}}&\geq \int_{\mathbb{R}^{n}}|\psi(s, y)|\, |y-x(s)|^{\omega}\,dy\geq  \int_{B_2}|\psi(s, y)| \,|y-x(s)|^{\omega}\,dy\\
&
\geq \frac{\psi(s,\overline{x}_{s})}{2}\int_{B_2} |y-x(s)|^{\omega}\,dy.\label{EstimationB2}
\end{align}
We just need to estimate the last integral following the cases given above. 
\begin{itemize}
\item Indeed, if $|\overline{x}_{s}-x(s)|>2\rho R$  then we have
$$\underset{y\in B_2\subset \mathcal{C}(R,\rho R)}{\min} |y- x(s)|^{\omega}\geq (\rho R)^{\omega}=\rho^{\omega}\,R^{\omega},$$
while if $|\overline{x}_{s}-x(s)|<R/2$, one has $
\displaystyle{\underset{y\in B_2\subset \mathcal{C}(R,\rho R)}{\min}|y-x(s)|^{\omega}\geq \frac{R^{\omega}}{2^\omega}}$. Applying these results to (\ref{EstimationB2}) we obtain
$$
\big((\zeta r)^{\alpha}+K s\big)^{(\omega-\gamma)/{\alpha}}\geq \frac{\psi(s,\overline{x}_{s})}{2} \,\rho^{\omega} \,R^{\omega}\,|B_2|
$$
and
$$
\big((\zeta r)^{\alpha}+K s\big)^{(\omega-\gamma)/{\alpha}}\geq \frac{\psi(s,\overline{x}_{s})}{2} \, \frac{R^{\omega}}{2^\omega}\,|B_2|,
$$
 and since $\rho>2$ we have the first desired estimate:
\begin{equation*}
\frac{C_1 \big((\zeta r)^{\alpha}+K s\big)^{(\omega-\gamma)/{\alpha}}}{\psi(s,\overline{x}_{s}) \,R^{\omega}} \geq |B_2|, \quad \mbox{with } C_1=2^{1+\omega}.
\end{equation*}

\item For the second case, since $R/2\leq |\overline{x}_{s}-x(s)|\leq 2\rho R$, we can write using the Cauchy--Schwarz inequality,
\begin{equation}\label{HolderInver}
\int_{B_2}|y-x(s)|^{\omega}\,dy\geq |B_2|^2\Big(\int_{B_2}|y-x(s)|^{-\omega}\,dy\Big)^{-1}.
\end{equation}
Now, observe that in this case we have $B_2\subset B(x(s), 5\rho R)$ and then
$$\int_{B_2}|y-x(s)|^{-\omega}\,dy\leq \int_{B(x(s), 5 \rho R)}|y-x(s)|^{-\omega}\,dy\leq v_n (5\rho R)^{n-\omega}.$$
Getting back to (\ref{HolderInver}) we have
$$
{\int_{B_2}}|y-x(s)|^{\omega}\,dy\geq |B_2|^2 \,v_n^{-1} \, (5 \rho R)^{-n+\omega},
$$
and we use this estimate in (\ref{EstimationB2}) to obtain
\begin{equation*}
\frac{C_2 \big((\zeta r)^{\alpha}+K s\big)^{(\omega-\gamma)/(2\alpha)} R^{n/2-\omega/2}}{\psi(s,\overline{x}_{s})^{1/2}}\geq |B_2|, \quad \mbox{where } C_2=(2\cdot  5^{n-\omega} v_n\,\rho^{n-\omega})^{1/2}.
\end{equation*}
\end{itemize}
Now, with estimates (\ref{FormEstima1}) and (\ref{FormEstima2}) at our disposal we can write
\begin{enumerate}
\item[$\bullet$] if $|\overline{x}_{s}-x(s)|>2\rho R$ or $|\overline{x}_{s}-x(s)|<R/2$ then
\begin{align*}
&\int_{\mathcal{C}(R,\rho R)} \frac{\psi(s, \overline{x}_s)-\psi(s, y)}{|\overline{x}_{s}-y|^{n+\alpha}}\,dy\\
&\qquad \geq  \frac{\psi(s, \overline{x}_s)}{2\,\rho^{n+\alpha}\, R^{n+\alpha}} \Big(v_n(\rho^n -1)R^n-\frac{C_1 \big((\zeta r)^{\alpha}+K s\big)^{(\omega-\gamma)/{\alpha}} } {\psi(s, \overline{x}_s)} \,R^{-\omega}\Big),
\end{align*}
\item[$\bullet$] if $R/2\leq |\overline{x}_{s}-x(s)|\leq 2\rho R$,
\begin{align*}
&\int_{\mathcal{C}(R,\rho R)}\frac{\psi(s, \overline{x}_s)-\psi(s, y)}{|\overline{x}_{s}-y|^{n+\alpha}}\,dy\\
&\qquad\geq  \frac{\psi(s, \overline{x}_s)}{2\,\rho^{n+\alpha}\, R^{n+\alpha}}\Big(v_n(\rho^n -1)R^n-\frac{C_2\big((\zeta r)^{\alpha}+K s\big)^{(\omega-\gamma)/(2\alpha)}R^{n/2-\omega/2}}{\psi(s, \overline{x}_s)^{1/2}}\Big).
\end{align*}
\end{enumerate}
If we set
$$
R= \big((\zeta r)^{\alpha}+K s\big)^{(\omega-\gamma)/(\alpha(n+\omega))}\, \psi(s, \overline{x}_s)^{-{1}/{(n+\omega)}},
 $$
 since we are working with small molecules we have $0<R\ll1$, and we obtain for all the previous cases the following estimate:
\begin{align*}
&\int_{\mathcal{C}(R,\rho R)}\frac{\psi(s, \overline{x}_s)-\psi(s, \overline{x}_s)}{|\overline{x}_{s}-y|^{n+\alpha}}\,dy\\
&\geq  \Big(\frac{v_n \,(\rho^n-1)-\sqrt{2v_n}\, (5\rho)^{(n-\omega)/{2}}}{2\,\rho^{n+\alpha}}\Big)\, \big((\zeta r)^{\alpha}+K s\big)^{-(\omega-\gamma)/(n+\omega)}\, \psi(s, \overline{x}_s)^{1+\alpha/(n+\omega)}.
\end{align*}
At this point, once the dimension $n$ and the parameters $\alpha, \omega$ are fixed, we obtain that the quantity
$$\mathfrak{C}=\frac{v_n \,(\textcolor{black}{\rho}^n-1)-\sqrt{2v_n}\, (5\textcolor{black}{\rho})^{(n-\omega)/2}}{2 \times {\textcolor{black}{\rho}}^{n+\alpha}},$$
can be made a small positive constant provided $\rho $ is large enough. Thus, and for all possible cases considered before, we have the following estimate for (\ref{Infty1}):
$$
\frac{d}{ds}\|\psi(s, \cdot)\|_{L^\infty}\leq -  \mathfrak{c}_1 \times  \mathfrak{C}\times \big((\zeta r)^{\alpha}+K s\big)^{-(\omega-\gamma)/(n+\omega)} \,\|\psi(s,\cdot)\|_{L^\infty}^{1+{\alpha}/(n+\omega)}.
$$
We recall now that the constant $K$ given in (\ref{ConditionK11}) can be small enough to write
$$
\frac{d}{ds}\|\psi(s, \cdot)\|_{L^\infty}\leq -K\, \Big(\frac{n+\gamma}{\alpha} \Big) \big((\zeta r)^{\alpha}+K s\big)^{-(\omega-\gamma)/(n+\omega)} \, \|\psi(s,\cdot)\|_{L^\infty}^{1+{\alpha}/(n+\omega)},$$
which is exactly formula (\ref{INEG_DIFF}).\\

The proof of the height condition is finished for regular molecules. In order to obtain the global result, remark that, for viscosity solutions we have $\Delta \psi(s_0, \overline{x})\leq 0$ at the points $\overline{x}$ where $\psi(s_0, \cdot)$ reaches its maximum value so we only need to study the term $\mathcal{L}^\alpha\psi(s_0, \overline{x})$ as it was done here.
We refer to \cite{Cordoba} and \cite{DCHSM} for more details. 
\begin{Remarque}
The constants obtained here do not depend on the molecule's size but only on the dimension $n$ and on parameters $\omega$,  $\gamma$ and $\alpha$.
\end{Remarque}
\subsubsection{$L^{1}$ estimate}\label{Secc_EstimationL1}
The $L^1$ control \eqref{SmallL1evolution} is now a direct consequence of an optimization over the parameter $D$ below (splitting threshold):
\begin{eqnarray*}
\int_{\mathbb{R}^n}|\psi(s_0,x)|dx&=&\int_{\{|x-x(s_0)|< D\}}|\psi(s_0,x)|dx+\int_{\{|x-x(s_0)|\geq D\}}|\psi(s_0,x)|dx\\
&\leq & v_n D^n \|\psi(s_0,\cdot)\|_{L^\infty}+D^{-\omega}\int_{\mathbb{R}^{n}}|\psi(s_0,x)||x-x(s_0)|^\omega dx.
\end{eqnarray*}
Now using the Concentration condition and the Height condition one has:
\begin{eqnarray*}
\int_{\mathbb{R}^n}|\psi(s_0,x)|dx&\leq & v_n \frac{D^n }{\left((\zeta r)^\alpha+Ks_0\right)^{\frac{n+\gamma}{\alpha}}}  +D^{-\omega}((\zeta r)^\alpha+Ks_0)^{\frac{\omega-\gamma}{\alpha}},
\end{eqnarray*}
where $v_n$ denotes the volume of the unit ball. An optimization over the real parameter $D$ yields:
\begin{equation*}
\|\psi(s_0,\cdot)\|_{L^1}\leq \frac{2v_n^{\frac{\omega}{n+\omega}}}{\big((\zeta r)^\alpha+Ks_0\big)^{\frac{\gamma}{\alpha}}}.
\end{equation*}
Theorem \ref{SmallGeneralisacion} is now completely proven. \hfill$\blacksquare$
\subsection{Iteration}\label{Secc_Iteration}
Once we have a good control over the quantities $\|\psi(s_0,\cdot)\|_{L^1}$ and $\|\psi(s_0,\cdot)\|_{L^\infty}$ \textcolor{black}{(from \eqref{SmallL1evolution} and \eqref{SmallLinftyevolution})}, by interpolation we obtain the following bound
$$\|\psi(s_0,\cdot)\|_{L^{p'}}\leq \|\psi(s_0,\cdot)\|_{L^1}^{\frac{1}{p'}} \|\psi(s_0,\cdot)\|_{L^\infty}^{1-\frac{1}{p'}}\leq C\left[\big((\zeta r)^\alpha+Ks_0\big)^{\frac{1}{\alpha}} \right]^{-n+\frac{n}{p'}-\gamma}.$$
We thus see with Theorem \ref{SmallGeneralisacion} that it is possible to control the $L^{p'}$ norm of the molecules $\psi$ from $0$ to a small time $s_0$, and \textcolor{black}{applying inductively} the same arguments we can extend the control from time $s_0$ to time $s_{1}$ with a small increment $s_{1}-s_{0}\sim \mathfrak{e}r^{\alpha}$.  Now we can see that the smallness of the parameters $\mathfrak{e}$, $r$ and of the time increments $s_0,s_1-s_0,...,s_N-s_{N-1}$ can be compensated by the number of iterations $N$: fix a small $0<r<1$ and iterate as explained before. Since each small time increment $s_0,s_1-s_0,...,s_N-s_{N-1} $ has order $\mathfrak{e} r^\alpha$, we have $s_{N}\sim N \mathfrak{e} r^\alpha$. Thus, we will stop the iterations as soon as $ N \mathfrak{e} r^\alpha\geq \mathcal{T}_0$. The number of iterations $N=N(\mathfrak{e},r)$ will depend on the smallness of the molecule's size $r$ and on the size of the initial data $\mu$ (controlled by the smallness of $\mathfrak{e}$), and it is enough to consider $N(\mathfrak{e},r)\sim \frac{\mathcal{T}_0}{\mathfrak{e} r^\alpha}$ in order to obtain this lower bound for $N(\mathfrak{e},r)$. Proceeding this way we will obtain $\|\psi(s_N,\cdot)\|_{L^{p'}}\leq C \mathcal{T}_0^{-n+\frac{n}{p'}-\gamma}<+\infty$, for all molecules of size $r$. Observe that once this estimate is available, for larger times it is enough to apply the maximum principle.

Finally, and for all $r>0$, we obtain after a time $T_0$ a $L^{p'}$ control for small molecules and we finish the proof of the Theorem \ref{TheoLp'control}. \hfill$\blacksquare$

\appendix
\mysection{Mild representation of weak solutions: proof of Lemma \ref{LEMME_MILD}}
\label{ANNEX_MILD}
We aim at proving that any weak solution $\theta $ to \eqref{EquationPrincipaleIntro} which belongs to $L^\infty_t(L^p_x)\cap L_t^p(\dot B_{p,x}^{\frac \alpha p,p}),\ p\in [2,+\infty[$ actually satisfies the Duhamel type representation \eqref{REP_MILD_W_SOL}. Fix $T>0$ and start from the identity:
\begin{equation}\label{WEAK_FORM}
\int_{\R^n} \theta_0(y)\varphi(0,y)dy=\int_0^T \int_{\R^n}  \theta(s,y) (-\partial_s+\mathcal L^\alpha)\varphi(s,y)dy ds
+\int_0^T  \int_{\R^n} \mathbb A_{[\theta]}\theta(s,y)\cdot \nabla \varphi(s,y)dyds,
\end{equation}
which holds for all $\varphi\in \mathcal{C}_0^\infty(]-T,T[\times \R^n) $. Similarly, for almost any $t\in ]0,T[$,
\begin{align}\label{WEAK_FORM_2}
\int_{\R^n} \theta_0(y)\varphi(0,y)dy=&\int_{\R^n}\theta(t,y) \varphi(t,y)dy+\int_0^t  \int_{\R^n}  \theta(s,y) (-\partial_s+\mathcal L^\alpha)\varphi(s,y)dy ds\notag\\
&+\int_0^t  \int_{\R^n} \mathbb A_{[\theta]}\theta(s,y)\cdot \nabla \varphi(s,y)dyds.
\end{align}
Observe that, for a fixed pair $(t,x)\in \R_+\times \R^n $, the mild representation \eqref{REP_MILD_W_SOL} would formally follow from \eqref{WEAK_FORM_2} taking therein $\varphi(s,y)$ as a smooth compactly supported approximation of $\mathfrak{p}_{t-s}^\alpha(y-x) $ recalling that 
\begin{equation}\label{HARMO_HK}
(-\partial_s+\mathcal L^\alpha)\mathfrak{p}_{t-s}^\alpha(y-x)=\delta_t(s)\delta_x(y)
\end{equation}
and passing to the limit. We actually claim that, in the current subcritical regime (recall that we have $\alpha=1+\varepsilon>0$) and thanks to the available controls (maximum principle, Besov energy estimates) this procedure can actually be justified and yields \eqref{REP_MILD_W_SOL} for almost any $(t,x)\in \R_+\times \R^n $. Precisely, define for a fixed $(t,x)\in \R_+^*\times \R^n$ and $\varepsilon>0$ the function $(s,y)\mapsto \varphi_{\varepsilon,(t,x)}(s,y)=\mathfrak{p}_{t-s+\varepsilon}^{\alpha,\varepsilon}(y-x) $ where for $\tau>0$, $\mathfrak p_{\tau}^{\alpha,\varepsilon}(z)=\Big(\mathfrak p_\tau^\alpha(\cdot) \mathds{1}_{|\cdot|\le \varepsilon^{-(2+\beta)}} \Big)\ast \phi_{\varepsilon^2}(z)$ for some $\beta>0 $ and where for $\eta>0 $, $ \phi_\eta$ is a standard mollifier. Then $\varphi_{\varepsilon,(t,x)}$ can be extended into a function in  $\mathcal{C}_0^\infty([-T,T[\times\R^n)$ and from \eqref{WEAK_FORM_2}:
\begin{align}\label{WEAK_FORM_3}
\int_{\R^n} \theta_0(y)\mathfrak p_{t}^\alpha(y-x)dy=&\theta(t,x)
+\int_0^t \int_{\R^n} \mathbb A_{[\theta]}\theta(s,y)\cdot \nabla \mathfrak p_{t-s}^\alpha(y-x)dy ds\notag\\
&+R_{\varepsilon}(t,x).
\end{align}
This will precisely give \eqref{REP_MILD_W_SOL} provided we prove that $R_{\varepsilon}(t,x) $ tends to 0 with $\varepsilon $ for almost any $(t,x)\in ]0,T]\times \R^n $. Write recalling \eqref{HARMO_HK}:
\begin{align*}
R_{\varepsilon}(t,x)=&\int_{\R^n}\theta(t,y) \varphi_{\varepsilon,(t,x)}(t,y)dy-\theta(t,x)+\int_0^t  \int_{\R^n}  \theta(s,y) (-\partial_s+\mathcal L^\alpha)\varphi_{\varepsilon,(t,x)}(s,y)dyds\\
&+\int_0^t  \int_{\R^n} \mathbb A_{[\theta]}\theta(s,y)\cdot \big( \nabla \varphi_{\varepsilon,(t,x)}(y)-\nabla \mathfrak p_{t-s}^\alpha(y-x)\big)dyds=:\sum_{i=1}^3 R_{\varepsilon,i}(t,x).
\end{align*}
Before going further we state a Lemma giving some quantitative controls for the stable like heat-kernel and its spatial sensitivities.
\begin{Lemme}[Controls for the stable-like heat kernel]\label{LEMME_HK}
Let \A{A} be in force and $\mathfrak p^\alpha $  denote the heat-kernel associated with ${\mathcal L}^\alpha $. There exists a function $\bar p $ satisfying that for all $z\in \R^n $, $\bar p(z)\le c(1+|z|)^{-(n+\alpha)} $ for some constant $C\ge 1$ (observe that $\bar p\in L^1(\R^n) $) and $C\ge 1$ such that for all $v>0, \bar z\in \R^n$, and any $\beta,\gamma \in \{0,1\},\ i\in \{1,\cdots,n\} $,
\begin{equation}\label{HK1}
v^{\frac n\alpha} |(\partial_v^\beta \partial_{x_i}^\gamma\mathfrak p_{v}^\alpha)(x)|_{x=v^{\frac 1\alpha}\bar z}\le Cv^{-(\beta+\frac\gamma\alpha)}\bar p(\bar z).
\end{equation}
Also, for all $\bar y\in \R^n, \bar z\in B(0,R) $, $\varepsilon>0 $:
\begin{equation}\label{SMALL_JUMPS}
v^{\frac n\alpha}\big|(\partial_v^\beta\partial_{x_i}^\gamma \mathfrak p_v^{\alpha})(x+\varepsilon^2\bar z) -(\partial_v^\beta\partial_{x_i}^\gamma\mathfrak p_v^{\alpha})(x)\big|_{x=v^{\frac 1\alpha}\bar y}\le C \varepsilon^2 v^{-(\beta+\frac \gamma \alpha+\frac1\alpha)} \bar p(\bar y).
\end{equation}
\end{Lemme}
The above results can be derived from the controls for the stable-like heat kernel and its sensitivities. In the isotropic stable case we can \emph{e.g.} refer to  \cite{meno:zhan:22}. The general case considered in Assumption \A{A}, where the tails of the L\'evy measure are upper-bounded by a stable one (which e.g. include the truncated stable kernel) could be handled following the approach introduced in \cite{wata:07}, consisting in precisely splitting the L\'evy measure into its small and large jumps part for the analysis of general stable heat kernels, and extended to the \textit{tempered} case in \cite{szto:10}.

Observe first that $R_{\varepsilon,1}(t,x)=\theta(t,\cdot)\ast \mathfrak p_\varepsilon^{\alpha,\varepsilon}(x) -\theta(t,x)$. It is then clear that $\|R_{\varepsilon,1}\|_{L^p([0,T],L^p(\R^n))} \underset{\varepsilon \rightarrow 0}{\longrightarrow} 0$.
Indeed,
\begin{align*}
\|R_{\varepsilon,1}\|_{L^p([0,T],L^p(\R^n))}\le \|\theta \ast \mathfrak p_\varepsilon^\alpha-\theta\|_{L^p([0,T],L^p(\R^n))}+\|\theta\|_{L^p([0,T],L^p(\R^n))} \|\mathfrak p_\varepsilon^{\alpha,\varepsilon}-\mathfrak p_\varepsilon^{\alpha}\|_{L^1(\R^n)}.
\end{align*}
The convergence to zero of the first term in the right-hand side follows from the fact that $\mathfrak p_\varepsilon^\alpha$ is an approximation of the identity whereas for the second one we have:
\begin{align*}
\|\mathfrak p_\varepsilon^{\alpha,\varepsilon}-\mathfrak p_\varepsilon^{\alpha}\|_{L^1}&= \int_{\R^n} |\mathfrak (\mathfrak p_\varepsilon^\alpha\mathds{1}_{|\cdot|\le \varepsilon^{-(2+\beta)}}-\mathfrak p_\varepsilon^\alpha)\ast \phi_{\varepsilon^2}(y)+(\mathfrak p_\varepsilon^\alpha\ast\phi_{\varepsilon^2}-\mathfrak p_\varepsilon^\alpha)(y)|dy\\
&\le \int_{|z|\ge \varepsilon^{-(2+\beta)}}\mathfrak p_\varepsilon^\alpha(z)dz+\int_{\R^n}|(\mathfrak p_\varepsilon^\alpha\ast\phi_{\varepsilon^2}-\mathfrak p_\varepsilon^\alpha)(y)|dy=:P_1^\varepsilon+P_2^\varepsilon.
\end{align*}
Thus, from \eqref{HK1} in Lemma \ref{LEMME_HK}, taking $v=\varepsilon, \beta=0 $ therein,
$$P_1^\varepsilon\le C\int_{|\bar z|\ge \varepsilon^{-(2+\beta+\frac 1\alpha)}} \bar p(\bar z) d\bar z\underset{\varepsilon \rightarrow 0}{\longrightarrow } 0. $$
Similarly, for $\phi_{\varepsilon^2}(\cdot)=\frac{1}{\varepsilon^{2n}} \phi(\frac\cdot{\varepsilon^2}) $ for a compactly supported in $B(0,R) $ (centered ball of radius $R$) non negative smooth $\phi $ such that $\displaystyle{\int_{\R^n} \phi(z)dz=\int_{B(0,R)} \phi(z)dz=1}$, we have
\begin{align*}
P_2^\varepsilon\le \int_{\R^n} \left|\int_{\R^n} \big(\mathfrak p_\varepsilon^{\alpha}(\varepsilon^{\frac 1\alpha}\bar y+z) -\mathfrak p_\varepsilon^{\alpha}(\varepsilon^{\frac 1\alpha}\bar y) \big)\phi_{\varepsilon^2}(z)  dz\right|\varepsilon^{\frac n\alpha} d\bar y\\
\le \int_{\R^n} \int_{\R^n}\big|\mathfrak p_\varepsilon^{\alpha}(\varepsilon^{\frac 1\alpha}\bar y+\varepsilon^2\bar z) -\mathfrak p_\varepsilon^{\alpha}(\varepsilon^{\frac 1\alpha}\bar y) \big| \phi(\bar z) d\bar z \varepsilon^{\frac n\alpha}d\bar y.
\end{align*}
Still from Lemma \ref{LEMME_HK}, equation \eqref{SMALL_JUMPS} with $v=\varepsilon , \beta=\gamma=0$ then gives
$P_2^\varepsilon \underset{\varepsilon \rightarrow 0}{\longrightarrow}0$ and therefore yields the stated convergence for $R_{\varepsilon,1} $.
Write now,
\begin{align*}
R_{\varepsilon,2}(t,x)&=\int_0^t  \theta(s,\cdot)\ast (-\partial_s+\mathcal L^\alpha)\mathfrak p_{t-s+\varepsilon}^{\alpha,\varepsilon}(x)ds=\int_0^t  \theta(s,\cdot)\ast (-\partial_s+\mathcal L^\alpha)[\mathfrak p_{t-s+\varepsilon}^{\alpha,\varepsilon}-\mathfrak p_{t-s+\varepsilon}^{\alpha}](x)ds,
\end{align*}
recalling \eqref{HARMO_HK} for the last equality. Then
\begin{align}\label{PREAL_BD_R2}
\|R_{\varepsilon,2}\|_{L^p([0,T],L^p(\R^n))}\le C\|\theta\|_{L^p([0,T],L^p(\R^n))} \int_0^T  \int_0^t   \|(-\partial_s+\mathcal L^\alpha)[\mathfrak p_{t-s+\varepsilon}^{\alpha,\varepsilon}-\mathfrak p_{t-s+\varepsilon}^{\alpha}]\|_{L^1(\R^n)}dsdt.
\end{align}
We can then reproduce the previous analysis and write:
\begin{align*}
&\|(\partial_s+\mathcal L^\alpha)\mathfrak p_{t-s+\varepsilon}^{\alpha,\varepsilon}-(\partial_s+\mathcal L^\alpha)\mathfrak p_{t-s+\varepsilon}^{\alpha}\|_{L^1(\R^n)}\\
= &\int_{\R^d} |(\partial_s+\mathcal L^\alpha)\mathfrak (\mathfrak p_{t+s-\varepsilon}^\alpha\mathds{1}_{|\cdot|\le \varepsilon^{-(2+\beta)}}-\mathfrak p_{t-s+\varepsilon}^\alpha)\ast \phi_{\varepsilon^2}(y)+(\partial_s+\mathcal L^\alpha)(\mathfrak p_{t-s+\varepsilon}^\alpha\ast\phi_{\varepsilon^2}-\mathfrak p_{t-s+\varepsilon}^\alpha)(y)|dy\\
\le &\int_{|z|\ge \varepsilon^{-(2+\beta)}}(|\partial_s \mathfrak p_{t-s+\varepsilon}^\alpha(z)|+| \mathfrak p_{t-s+\varepsilon}^\alpha(z)| \|\mathcal L^\alpha \phi_{\varepsilon^2}\|_{L^1})dz+\int_{\R^n}|(\partial_s+\mathcal L^\alpha)(\mathfrak p_{t-s+\varepsilon}^\alpha\ast\phi_{\varepsilon^2}-\mathfrak p_{t-s+\varepsilon}^\alpha)(y)|dy\\
=:&Q_1^\varepsilon+Q_2^\varepsilon.
\end{align*}
Let us first consider here $Q_1^\varepsilon$. Equation \eqref{HK1} with $v=t-s+\varepsilon $ and $\beta=1,\gamma=0 $ yields
\begin{equation}\label{CTR_DER_TEMP_DENS}
(t-s+\varepsilon)^{\frac n\alpha} |\partial_s \mathfrak p_{t-s+\varepsilon}^\alpha(\bar x)|_{\bar x=(t-s+\varepsilon)^{\frac 1\alpha}\bar z}\le C(t-s+\varepsilon)^{-1}\bar p(\bar z).
\end{equation}
By homogeneity, we also get $\|\mathcal L^\alpha \phi_{\varepsilon^2}\|_{L^1}\le C\varepsilon^{-2\alpha}$. Hence,
\begin{align*}
Q_1^\varepsilon\le C\int_{ |\bar z|\ge (t-s+\varepsilon)^{-\frac1\alpha }\varepsilon^{-(2+\beta)}}((t-s+\varepsilon)^{-1}+\varepsilon^{-2\alpha} )\bar p(\bar z)d\bar z.
\end{align*}
Recalling now that $\bar p(\bar z)\le C(1+|\bar z|)^{-(n+\alpha)} $ we get:
\begin{align*}
Q_1^\varepsilon\le C(\varepsilon^{\alpha(2+\beta)}+(t-s+\varepsilon)\varepsilon^{\alpha \beta})\underset{\varepsilon \rightarrow 0}{\longrightarrow}0.
\end{align*}
Write now for $Q_2^\varepsilon$,
\begin{align*}
Q_2^\varepsilon&\le \int_{\R^n} \left|\int_{\R^n} (\partial_s+\mathcal L^\alpha)\big(\mathfrak p_{t-s+\varepsilon}^{\alpha}(\bar x+z) -\mathfrak p_{t-s+\varepsilon}^{\alpha}(\bar x )\big)\phi_{\varepsilon^2}(z)  dz\right|_{\bar x=(t-s+\varepsilon)^{\frac 1\alpha}\bar y}(t-s+\varepsilon)^{\frac n\alpha} d\bar y\\
&\le \int_{\R^n} \int_{\R^n}\big|(\partial_s +\mathcal L^\alpha)(\mathfrak p_{t-s+\varepsilon}^{\alpha}(\bar x+\varepsilon^2\bar z) -\mathfrak p_{t-s+\varepsilon}^{\alpha}(\bar x) \big|_{\bar x=(t-s+\varepsilon)^{\frac 1\alpha}\bar y} \phi(\bar z) d\bar z (t-s+\varepsilon)^{\frac n\alpha}d\bar y.
\end{align*}
Similarly to the previous controls on the heat kernel we derive from \eqref{SMALL_JUMPS} with $v=t-s+\varepsilon,\beta=1,\gamma=0 $ (recall that $\bar z\in B(0,R) $) and \eqref{HARMO_HK}:
\begin{equation}\label{SENSI_TIME_DER}
(t-s+\varepsilon)^{\frac n\alpha}|(\partial_s +\mathcal L^\alpha)(\mathfrak p_{t-s+\varepsilon}^{\alpha}(\bar x+\varepsilon^2\bar z) -\mathfrak p_{t-s+\varepsilon}^{\alpha}(\bar x)) \big|_{\bar x=(t-s+\varepsilon)^{\frac 1\alpha}\bar y}\le C (t-s+\varepsilon)^{-(1+\frac 1\alpha)}\bar p(\bar y)\varepsilon^2,
\end{equation}
which again readily gives $ Q_2^\varepsilon \le C \varepsilon^{1-\frac 1\alpha}$. Plugging the above controls into \eqref{PREAL_BD_R2} we derive the required convergence for the term $R_{\varepsilon,2} $. Let us now turn to $R_{\varepsilon,3} $.
\begin{align*}
\|R_{\varepsilon,3}\|_{L^1([0,T],L^1(\R^n))}\le& \left\|\int_0^\cdot  \int_{\R^n} \mathbb A_{[\theta]}\theta(s,y)\cdot \big( \nabla \varphi_{\varepsilon,(t,\cdot)}(y)-\nabla \mathfrak p_{t-s}^\alpha(y-\cdot)\big)dyds\right\|_{L^1([0,T],L^1(\R^n))}\\
\le& \int_0^T  \int_0^t  \|\mathbb A_{[\theta]}\theta(s,\cdot)\|_{L^1}\|\nabla\mathfrak p_{t-s+\varepsilon}^{\alpha,\varepsilon}-\nabla \mathfrak p_{t-s+\varepsilon}^{\alpha}\|_{L^1}dsdt\\
\le &\|\theta\|_{L_t^\infty (L^2_x)}\|\mathbb A_{[\theta]}\|_{L_t^\infty (L_x^2)}\int_0^T  \int_0^t  \|\nabla\mathfrak p_{t-s+\varepsilon}^{\alpha,\varepsilon}-\nabla \mathfrak p_{t-s+\varepsilon}^{\alpha}\|_{L^1}dsdt\\
\le &C\|\theta_0\|_{L^2}^2\int_0^T  \int_0^t \|\nabla\mathfrak p_{t-s+\varepsilon}^{\alpha,\varepsilon}-\nabla \mathfrak p_{t-s+\varepsilon}^{\alpha}\|_{L^1}ds dt,
\end{align*}
using (\ref{HypoLPBorne1}) and by  maximum principle (\ref{MaximumPrinciple1}) for the last inequality. The contribution of the term $\|\nabla\mathfrak p_{t-s+\varepsilon}^{\alpha,\varepsilon}-\nabla \mathfrak p_{t-s+\varepsilon}^{\alpha}\|_{L^1}$ can be handled just as the corresponding one for $R_{2,\varepsilon} $, \emph{i.e.} the quantity $\|(-\partial_s+\mathcal L^\alpha)[\mathfrak p_{t-s+\varepsilon}^{\alpha,\varepsilon}-\mathfrak p_{t-s+\varepsilon}^{\alpha}]\|_{L^1(\R^n)}$ in \eqref{PREAL_BD_R2}, noting that the associated time singularity will be milder. Namely, its order will be $(t-s+\varepsilon)^{-\frac 1\alpha} $ (usual time singularity associated with the spatial gradient of the stable heat-kernel) instead of $(t-s+\varepsilon)^{-1} $ (usual time singularity for the time derivative or the fractional operator applied to the heat kernel) for the previous case. Applying \eqref{SMALL_JUMPS} with $ \beta=0,\gamma=1$, the associated controls then write
\begin{align*}
\label{CTR_DER_SPACE_DENS}
(t-s+\varepsilon)^{\frac n\alpha} |\nabla \mathfrak p_{t-s+\varepsilon}^\alpha(\bar x)|_{\bar x=(t-s+\varepsilon)^{\frac 1\alpha}\bar z}\le C(t-s+\varepsilon)^{-\frac 1\alpha}\bar p(\bar z),\\
(t-s+\varepsilon)^{\frac n\alpha}|\nabla(\mathfrak p_{t-s+\varepsilon}^{\alpha}(\bar x+\varepsilon^2\bar z) -\mathfrak p_{t-s+\varepsilon}^{\alpha}(\bar x) \big|_{\bar x=(t-s+\varepsilon)^{\frac 1\alpha}\bar y}\le C (t-s+\varepsilon)^{-\frac 2\alpha}\bar p(\bar y)\varepsilon^2,
\end{align*}
instead of \eqref{CTR_DER_TEMP_DENS}, \eqref{SENSI_TIME_DER}.\\
This gives the convergence to 0 of $R_{\varepsilon,3}$ and gives the representation \eqref{REP_MILD_W_SOL} for almost all $(t,x)\in [0,T]\times \R^d $. In the current sub-critical regime, the existence of a continuous in time version is immediate. This therefore  completes the proof of Lemma \ref{LEMME_MILD}. \hfill $\blacksquare$


\quad\\

\begin{multicols}{2}
\begin{minipage}[r]{90mm}
Diego \textsc{Chamorro}\\[3mm]
{\footnotesize
Laboratoire de Mod\'elisation\\ Math\'ematique d'Evry\\ 
(LaMME), UMR CNRS 8071\\ 
Universit\'e d'Evry Val d'Essonne\\[2mm]
23 Boulevard de France\\
91037 Evry Cedex\\[2mm]
diego.chamorro@univ-evry.fr
}
\end{minipage}
\begin{minipage}[r]{80mm}
St\'ephane \textsc{Menozzi}\\[3mm]
{\footnotesize
Laboratoire de Mod\'elisation\\ Math\'ematique d'Evry\\ 
(LaMME), UMR CNRS 8071\\
Université d'Evry Val d'Essonne\\[2mm]
23 Boulevard de France\\
91037 Evry Cedex\\and \\
Laboratory of Stochastic Analysis, HSE\\ 
Moscow, Russia\\[2mm]
stephane.menozzi@univ-evry.fr
}
\end{minipage}
\end{multicols}

\end{document}